\magnification 1200
\input amssym.def
\input amssym.tex
\parindent = 40 pt
\parskip = 12 pt
\font \heading = cmbx10 at 12 true pt
 at 22 true pt
 at 19 true pt
 at 7 true pt
\def \R{{\bf R}}

\centerline{\heading The Asymptotic Behavior of Degenerate Oscillatory}
\centerline{\heading Integrals in Two Dimensions}
\rm
\line{}
\line{}
\centerline{\heading Michael Greenblatt}
\line{}
\centerline{February 3, 2009}
\baselineskip = 12 pt
\font \heading = cmbx10 at 14 true pt
\line{}
\line{}
\noindent{\bf 1. Introduction}

\vfootnote{}{This research was supported in part by NSF grant DMS-0654073}In this paper we 
are interested in the following type of oscillatory integral. Suppose $S(x,y)$ is a 
smooth real-valued function defined in a neighborhood of the origin in $\R^2$, and $\phi(x,y) \in
C_c^{\infty} (\R^2)$ is real-valued and supported in a small neighborhood of the origin. We 
define
$$J_{S, \phi}(\lambda) = \int_{\R^2} e^{i \lambda S(x,y)} \phi(x,y)\,dx\,dy\eqno (1.1)$$
Here $\lambda$ is a real parameter and we want to understand the behavior of $J_{S, \phi}(\lambda)$
as $\lambda \rightarrow +\infty$. Oscillatory integrals of the form $(1.1)$ and their higher-dimensional
analogues come up frequently in several areas of analysis, including PDE's, mathematical
physics, and harmonic analysis. For example, such oscillatory integrals arise when analyzing the decay 
of Fourier transforms of surface-supported measures such as in [IKM].
We refer to chapter 8 of [S] for an overview of such issues. The stability of 
oscillatory integrals $(1.1)$ under perturbations of the phase function $S(x,y)$ is connected to 
various issues in complex geometry and has been considered for example in [PSSt] and [V].

Since one can always factor out an $e^{i \lambda S(0,0)}$, it does
no harm to assume that $S(0,0) = 0$. If $\nabla S(0,0) \neq 0$, in a small enough neighborhood of the 
origin one can integrate by parts arbitrarily many times in $(1.1)$ and get that $J_{S, \phi}(\lambda)$
decays faster than $C_N\lambda^{-N}$ for any $N$. Hence in this paper we always assume that the origin 
is a critical point for $S$; that is, we assume that
$$S(0,0) = 0,\,\,\,\,\,\,\,\nabla S(0,0) = 0$$
In this case, if $S(x,y)$ is real-analytic, using resolution 
of singularities (see [G1] for an elementary proof) one always has an asymptotic expansion
$$J_{S, \phi}(\lambda) \sim \sum_{j = 0}^{\infty}(d_j(\phi)\lambda^
{-s_j}  + d_j'(\phi)\ln(\lambda) \lambda^{-s_j}) \eqno (1.2)$$
Here $\{s_j\}$ is an increasing arithmetic progressions of positive rational numbers 
independent of $\phi$ deriving from the resolution of singularities of $S$. We always assume $s_0$ 
is chosen to be minimal such that in any sufficiently small neighborhood $U$ of the origin 
$d_0(\phi)$ or $d_0'(\phi)$ is nonzero for some $\phi$ supported in $U$. In this paper, we will
give explicit formulas for the leading term of $(1.2)$ once one is in certain coordinate
systems which we call ``superadapted", in analogy with the adapted coordinate systems of [V].
In the smooth case, we will find appropriate weaker analogues.

\noindent {\bf Definition 1.1.}  The {\it oscillatory index} of $S$ is defined to be $s_0$.
If in any small neighborhood of the origin there is some $\phi$ for which $d_0'(\phi)$ is nonzero, 
then we say $s_0$ has multiplicity 1. Otherwise, we say it has multiplicity zero. 

In the case where $S$ is a smooth function whose 
Hessian determinant at the origin is nonvanishing, one can do a smooth
coordinate change such that $S(x,y)$ becomes $\pm x^2 \pm y^2$ in the new coordinates. Then one
can use the well-known one dimensional theory (see Chapter 8 of [S]) and explicitly obtain an asymptotic 
expansion for $J_{S, \phi}(\lambda)$; the leading term will be given by ${2 \pi i \over 
D^{{1 \over 2}}} \phi(0,0)\lambda ^{-1}$ where $D$ denotes the Hessian determinant of $S$ at the 
origin. Hence our concern in this paper will be when $D = 0$; that is, when $S$ has a 
degenerate critical point at the origin. 

In the real-analytic situation, there is a close relationship between $J_{S,\phi}(\lambda)$ and the 
function $I_{S,\phi}(\epsilon)$ defined by
$$I_{S, \phi}(\epsilon) =  \int_{\{(x,y): 0 < S(x,y) < \epsilon\}} \phi(x,y) \,dx\,dy \eqno (1.3)$$
Analogous to $(1.2)$, when $S(0,0) = 0$ the functions $I_{S, \phi}(\epsilon)$ and $I_{-S, \phi}
(\epsilon)$ have asymptotic expansions which we may write as 
$$I_{S, \phi}(\epsilon) \sim \sum_{j = 0}^{\infty} (c_j(\phi)\epsilon^{r_j} + c_j'(\phi)\ln(\epsilon)
\epsilon^{r_j}) \eqno (1.4a)$$
$$I_{-S, \phi}(\epsilon) \sim \sum_{j = 0}^{\infty} (C_j(\phi)\epsilon^{r_j} + C_j'(\phi)\ln(\epsilon)
\epsilon^{r_j}) \eqno (1.4b)$$

Analogous to before, $\{r_j\}$ is an increasing arithmetic progression of positive 
rational numbers independent of $\phi$ deriving from the resolution of singularities of $S$, and
$r_0$ is chosen to be minimal such that in any sufficiently small neighborhood $U$ of the origin 
there is some $\phi$ supported in $U$ for which at least one of $c_0(\phi)$, 
$c_0'(\phi)$, $C_0(\phi)$, or $C_0'(\phi)$ is nonzero.

Using well-known methods (see Ch 7 of [AGV]), when the oscillatory index is less than 1, corresponding 
to the degenerate case, $r_0 = s_0$ and the coefficient of the leading term of $(1.2)$ can always be 
expressed in terms of those of the corresponding terms of $(1.4a)$ and 
$(1.4b)$. Often $I_{S, \phi}$ and $I_{-S, \phi}$ are easier to deal with than $J_{S, \phi}$ due to the 
absence of cancellations which can make it difficult to find lower bounds for $J_{S, \phi}$ directly. 

\noindent The $I_{S, \phi}$ are closely related to the $H_{S,\,U}$ defined by
$$H_{S,\,U}(\epsilon) = |\{x \in U: 0 < S(x,y) < \epsilon\}| \eqno (1.5)$$
Here $U$ is a small open set containing the origin, and the goal is to understand how $H_{S,\, U}
(\epsilon)$ behaves as $\epsilon \rightarrow 0$. Many of our results concerning the $I_{S, \phi}
(\epsilon)$ will immediately imply corresponding results about the $H_{S,\,U}(\epsilon)$. Specifically, 
one chooses $\phi_1$ supported in $U$ and equal to 1 outside some $\delta$-neighborhood of the boundary
of $U$, then
chooses $\phi_2$ equal to 1 on $U$ and supported on a $\delta$-neighborhood of $U$. One compares the
theorems for $I_{S,\phi_1}(\epsilon)$ and $I_{S,\phi_2}(\epsilon)$ and then lets $\delta$ go to zero if
necessary.
In the case where $S$ has an isolated zero at
the origin, for $\epsilon$ small enough the set $\{x \in U: 0 < S(x,y) < \epsilon\}$ will be a subset of
a set where the theorems hold, so one can simply take $\phi = 1$ and then $H_{S,\,U}(\epsilon) 
= I_{S,\phi}(\epsilon)$ for such $\epsilon$. 

In [V], Varchenko developed some ideas that went a long way towards understanding the case of degenerate
real-analytic phase. To describe his work, we need some pertinent definitions.

\noindent {\bf Definition 1.2.} Let $S(x,y) = \sum_{a,b} s_{ab}x^ay^b$ denote the 
Taylor expansion of $S(x,y)$ at the origin. Assume there is at least one $(a,b)$ for which $s_{ab}$ is
nonzero. For any $(a,b)$ for which $s_{ab} \neq 0$, let $Q_{ab}$ be the quadrant $\{(x,y) \in \R^2: 
x \geq a, y \geq b \}$. Then the {\it Newton polygon} $N(S)$ of $S(x,y)$ is defined to be 
the convex hull of the union of all $Q_{ab}$.  

In general, a Newton polygon consists of finitely many (possibly zero) bounded edges of negative slope
as well as an unbounded vertical ray and an unbounded horizontal ray.

\noindent {\bf Definition 1.3.} The {\it Newton distance} $d(S)$ of $S(x,y)$ is defined to be 
$\inf \{t: (t,t) \in N(S)\}$.

Throughout this paper, we will use the $(t_1,t_2)$ coordinates to write equations of lines relating
to Newton polygons, so as to distinguish from the $x$-$y$ variables of the domain of $S(x,y)$. The line 
in the $t_1$-$t_2$ plane with equation $t_1 = t_2$ comes up so frequently it has its own name:

\noindent {\bf Definition 1.4:} The {\it bisectrix} is the line in the $t_1$-$t_2$ plane with equation
$t_1 = t_2$.

A key role in the above theorems as well as our theorems to follow is played by the following polynomials.

\noindent {\bf Definition 1.5}. Suppose $e$ is a compact edge of $N(S)$. Define $S_e(x,y)$
by $S_e(x,y) = \sum_{(a,b) \in e} s_{ab} x^{a}y^{b}$. In other words $S_e(x,y)$ is the sum of the terms
of the Taylor expansion of $S$ corresponding to $(a,b) \in e$. If $S(x,y)$ is real-analytic, we 
use the same terminology when $e$ is the vertical or horizontal ray of $N(S)$. 

In [V], Varchenko showed that when $S$ is real-analytic the oscillatory index is always at most 
${1 \over d(S)}$, and that there is necessarily a coordinate system in 
which it is actually equal to ${1 \over d(S)}$. He also showed that the coordinate change to such
coordinates can 
always be made of the form $(x,y) \rightarrow (x, y - f(x))$ or $(x,y) \rightarrow (x - f(y), y)$ 
for real analytic $f$. Coordinate systems where $d(S)$ achieves the maximum possible value are referred 
to as ``adapted coordinates". He also showed that the multiplicity of the oscillatory index is equal 
to 1 if and only if there are adapted coordinates where the bisectrix intersects $N(S)$ at
a vertex. Otherwise the multiplicity is 0; the leading term of $(1.2)$ will not have the  
$\ln(\lambda)$ factor in it. The issue of finding an expression for the leading coefficient $d_0(\phi)$ or 
$d_0'(\phi)$ is not treated in [V]. However, in the case where the $S(x,y)$ has a critical point of
finite Milnor number at the origin, it is shown in [V] that the leading coefficient $d_0(\phi)$ or
$d_0'(\phi)$ is some
fixed multiple of $\phi(0,0)$ depending on the phase; precisely which multiple is not determined. 
These results were later extended in [Sh] to convex finite-type functions. 

Smooth analogues are proven in [IM] and [IKM]. In [IM] it is shown that adapted coordinates exist in 
the smooth case. In [IKM] it is shown that in smooth adapted coordinates
one has the estimate $|J_{S, \phi}(\lambda)| < C\ln|\lambda| |\lambda|^{-{ 1 \over d(S)}}$ for large
$|\lambda|$. There are also operator versions of
such results. For example, in [PS] it is proven that $||\int e^{i\lambda S(x,y)}\phi(x,y) 
f(y)\,dy||_{L^2} $ $< C|\lambda|^{-{1 \over 2 d(\tilde{S})}}||f||_{L^2}$ for real analytic phase, where 
$\tilde{S}(x,y)= S(x,y) - S(0,y) - S(x,0)$. The exponent ${1 \over 2d(\tilde{S})}$ is sharp. 
Generalizations to smooth phase were proven in [R] and [G2]. These results use subdivisions
into curved regions as will be done here. However, there are significant differences since one
gets stronger results for operators; in particular, adapted coordinates are not needed. 

Our theorems below will require us to be in certain adapted coordinate systems which we call
``superadapted" coordinate systems:

\noindent {\bf Definition 1.6}. One is in {\it superadapted coordinates} if whenever $e$ is a compact
edge of $N(S)$ intersecting the bisectrix, both of the functions $S_e(1,y)$ and $S_e(-1,y)$ have 
no real zero of order $d(S)$ or greater other than possibly $y = 0$.

It can be shown that an equivalent definition is obtained by stipulating the same condition on $S_e(x,1)$ 
and $S_e(x,-1)$ instead of $S_e(1,y)$ and $S_e(-1,y)$; we choose the $y$-variable for definiteness. In
section 7 we will prove any phase function can be put in superadapted coordinates using some ideas 
from two-dimensional resolution of singularities. 

\noindent {\bf Lemma 1.0.} In a superadapted coordinate system, a critical point of $S(x,y)$ at the 
origin is nondegenerate if and only if $d(S) = 1$.

\noindent {\bf Proof.} Write $S(x,y) = ax^2 + bxy + cy^2 + O(|x|^3 + |y|^3)$. The only way
$d(S)$ could be greater than 1 is for either $a$ and $b$ to both be zero, or for $c$ and $b$ to both
be zero. In either case, the Hessian at the origin is zero and the phase is degenerate. So we assume
that $d(S) = 1$, and we will show that in a superadapted coordinate system the phase is nondegenerate.

First consider the case where $N(S)$ has a vertex at $(1,1)$. Then $b \neq 0$ and either $a$ or $c$ 
is zero. Suppose $a = 0$ but $c \neq 0$; we claim that this
implies the coordinate system is not superadapted. For in this case there is an edge $e$ connecting
$(1,1)$ and $(0,2)$. Then $S_e(x,y) = bxy + cy^2$ and thus $S_e(1,y) = by + cy^2$ has a zero at 
$-{b \over c} \neq 0$, inconsistent with the definition of superadapted. Thus in a superadapted
coordinate system, if $N(S)$ has a vertex at $(1,1)$ and $a = 0$, then $c = 0$ and thus the 
Hessian is nonzero at the origin; the phase is nondegenerate. The case where $c = 0$ but $a \neq 0$
leads to a similar contradiction. We conclude that if $N(S)$ has a vertex at $(1,1)$ and $d(S) = 1$
then the phase is nondegenerate.

Next, consider the case where $(1,1)$ is in the interior of an edge $e$ of $N(S)$. In this case, the
endpoints of $e$ are $(2,0)$ and $(0,2)$. Hence $a$ and $c$ are nonzero. In a superadapted coordinate
system, one must have that $S(1,y) = a + by + cy^2$ has no real zeroes other than $y = 0$. Since 
$a \neq 0$, this is equivalent to $a + by + cy^2$ having no real zeroes at all, which happens 
exactly when $b^2 < 4ac$. This is equivalent to the Hessian determinant at the origin being nonzero,
and thus the phase is nondegenerate in this situation too. This completes the proof of Lemma 1.0.

We now come to our theorems. We will use the shorthand $d$ to denote the Newton distance 
$d(S)$. If $N(S)$ intersects the bisectrix in the interior of
an edge, bounded or unbounded, we denote this edge by $e_0$ and its slope by $-{1 \over m}$, where
$0 \leq m \leq \infty$. We use the shorthand $S_0(x,y)$ to denote $S_{e_0}(x,y)$. 

Theorem 1.1 is our main result. It gives explicit expressions for the leading term of $(1.4a)$ in 
superadapted coordinates. Applying the theorem to $-S$ gives analogous 
formulas for the expansion $(1.4b)$. As indicated above, $(1.4a)-(1.4b)$ directly
imply formulas for the leading term of the asymptotic expansion $(1.2)$ in the degenerate 
case; these are given in Theorem 1.2. 

\noindent {\bf Theorem 1.1.} Suppose $S(x,y)$ is a smooth
phase function in superadapted coordinates with $d > 1$. If the function $\phi(x,y)$ is supported in 
a sufficiently small neighborhood of the origin, then $r_0 = {1 \over d}$ and the following hold.

\noindent {\bf a)} Suppose the bisectrix intersects $N(S)$ in the interior of a compact 
edge. Define the function $S_0^+(x,y)^{-{1 \over d}}$ to be $S_0(x,y)^{-{1 \over d}}$ when $S_0(x,y) > 0$ and
zero otherwise. Then  we have
$$\lim_{\epsilon \rightarrow 0} {I_{S,\phi}(\epsilon) \over \epsilon^{1 \over d}} = (m + 1)^{-1}
\phi(0,0) \int_{-\infty}^{\infty} (S_0^+(1,y)^{-{ 1 \over d}} + S_0^+(-1,y)^{-{ 1 \over d}})\,dy 
\eqno (1.6)$$
In particular, if $S(x,y)$ is real-analytic then the coefficient $c_0'(\phi)$ in $(1.4a)$ is always
zero and $c_0(\phi)$ is given by $(1.6)$.

\noindent {\bf b)} Suppose the bisectrix intersects $N(S)$ at a 
vertex $(d,d)$. Let $s_{dd}x^dy^d$ denote the corresponding term of the Taylor 
expansion of $S$; hence $s_{dd} \neq 0$. Denote the slopes of the two edges of $N(S)$ meeting at 
$(d,d)$ by $s_1$ and $s_2$, where $-\infty \leq s_2 <  s_1 \leq 0$. Then 
$$\lim_{\epsilon \rightarrow 0} {I_{S,\phi}(\epsilon) \over \epsilon^{1 \over d}\ln(\epsilon)} = 
\eta(S)|s_{dd}|^{-{1 \over d}}\phi(0,0)({1 \over s_1 - 1} - {1 \over s_2 - 1}) \eqno (1.7)$$
Here $\eta(S)= 4$ if $s_{dd} > 0$ and $d$ is even, $\eta(S) = 2$ if $d$ is odd, and $\eta(S) = 0$ if 
$s_{dd} < 0$ and $d$ is even. 

\noindent {\bf c)} Suppose $S(x,y)$ is real-analytic and the bisectrix intersects $N(S)$ in the 
interior of the horizontal ray.  Write $S_{0}(x,y) = a(x)y^d$ where $a(x)$ is real-analytic. Let 
$\alpha(x)$ denote the one-dimensional measure of $\{y: 0 \leq a(x)y^d \leq 1\}$. In particular,
$\alpha(x) = |a(x)|^{-{1 \over d}}$ when $d$ is odd. Then $c_0'(\phi) = 0$, and $c_0(\phi)$ is given by
$$c_0(\phi) = \int_{-\infty}^{\infty} \alpha(x)\phi(x,0)\,dx \eqno (1.8)$$
The case where the bisectrix intersects the interior of the vertical ray has the analogous formula.

As can be seen, parts a), b) and c) give quite different formulas. Correspondingly, in our
subsequent theorems we break up into three cases. Case
1 is when the bisectrix intersects $N(S)$ in the interior of a compact edge, case 2 is when the 
bisectrix 
intersects $N(S)$ at a vertex $(d,d)$, and case 3 is when the bisectrix intersects $N(S)$ in the
interior of one of the unbounded rays. Notice that in cases 1 and 3, for a given $\phi(x,y)$ the 
expressions of Theorem 1.1 depend only on
$S_0(x,y)$, and that in case 2 it depends on $s_{dd}x^dy^d$ as well as the slopes of the edges of $N(S)$
meeting at $(d,d)$. 

Our next theorem gives the oscillatory integral version of Theorem 1.1 for the real-analytic case. 

\noindent {\bf Theorem 1.2.} Assume $S(x,y)$ is real-analytic and is in superadapted coordinates 
with $d > 1$. Then $s_0 = {1 \over d}$. In case 1 and
3, the coefficient $d_0'(\phi)$ of $(1.2)$ is always zero and $d_0(\phi)$ is given by
$$ d_0(\phi) = {\Gamma({1 \over d}) \over d}(e^{i {\pi \over 2d}} c_0(\phi) + e^{-i{ \pi \over 2d}}
C_0(\phi))\eqno (1.9a)$$
In case 2, one has
$$ d_0'(\phi) = -{\Gamma({1 \over d}) \over d}(e^{i {\pi \over 2d}} c_0'(\phi) + e^{-i{ \pi \over 2d}}
C_0'(\phi)) \eqno (1.9b)$$

\noindent {\bf Proof.} We will be sketchy here since the method for proving Theorem 1.2 from Theorem 1.1
is well-known; we refer to chapter 7 of [AGV] for more details. Note that 
$$J_{S,\phi}(\lambda) = \int_0^{\infty}({\partial_{\epsilon}}I_{S,\phi} (\epsilon))e^{i\lambda \epsilon}
\,d\epsilon + \int_0^{\infty} ({\partial_{\epsilon}}I_{-S,\phi}(\epsilon)) e^{-i\lambda \epsilon} 
\,d\epsilon \eqno (1.10a)$$
By [F], if $\gamma \in C_c(\R)$ with $\gamma(t) = 1$ near 0 and if $\alpha > -1$, for any $l$ we have 
$$\int_0^{\infty} e^{i \lambda t}t^{\alpha}\ln(t)^m \gamma(t)\,dt = {\partial^m \over \partial \alpha^m} 
{\Gamma(\alpha+1) \over (-i\lambda)^{\alpha + 1}} + O(\lambda^{-l})\eqno (1.10b)$$
Inserting $(1.4a)$ and $(1.4b)$ into $(1.10a)$ and using Theorem 1.1 and (1.10b) gives the theorem.

\noindent {\bf Comment 1.} One does need that $d > 1$ for Theorem 1.2 to hold. When $S(x,y)$ has a nondegenerate
saddle critical point, there are coordinates where $S(x,y) = xy$. This falls under case 2, and
one has that $c_0'(\phi) = C_0'(\phi)$. This means that the two terms of $(1.9b)$ will cancel. And
in fact when $\phi(0,0) \neq 0$, $|I_{S,\phi}(\epsilon)|,\,\,|I_{-S,\phi}(\epsilon)| \sim |\epsilon
\ln(\epsilon)|$, while $|J_{S,\phi}(\lambda)| \sim \lambda^{-1}$. 
 
\noindent {\bf Comment 2.} In cases 1 and 2, the expressions $(1.9a)$ and $(1.9b)$ for $d_0(\phi)$
and $d_0'(\phi)$ will always be nonzero when $d > 1$ and $\phi(0,0) \neq 0$. This is because the 
expressions given by Theorem 1.1 for $c_0(\phi)$, 
$C_0(\phi)$, $c_0'(\phi)$, and $C_0'(\phi)$ are real multiples of $\phi(0,0)$, while the ratio of the 
$e^{i {\pi \over 2d}}$ and $e^{-i {\pi \over 2d}}$ factors is never real when $d > 1$. 

\noindent {\bf Comment 3.} In any dimension, when the phase satisfies an appropriate 
nondegeneracy condition there are reasonably explicit formulas for the leading coefficient of the 
leading term of the asymptotic expansion of oscillatory integrals such as $(1.1)$. Such formulas are 
proven in [DS] and [DNS].  

\noindent Next, we give some less precise $C^{\infty}$ analogues for the $I_{S,\phi}(\epsilon)$.
The lower bounds involve $I_{|S|,\phi}(\epsilon) = I_{S,\phi}(\epsilon) + I_{-S,\phi}(\epsilon) = 
\int_{\{(x,y): |S(x,y)| < \epsilon\}} \phi(x,y)\,dx\,dy$. 

\noindent {\bf Theorem 1.3a.} Suppose now that $S(x,y)$ is a smooth phase function in superadapted
coordinates with $d > 1$. If $\phi$ is supported in a sufficiently small neighborhood of the
origin, then there is a positive $ B_{S,\phi}$ such that:

\noindent In cases 1 and 3 one has 
$$I_{S, \phi}(\epsilon) < B_{S,\phi}\epsilon^{1 \over d} \eqno (1.11)$$
In case 2 one has 
$$I_{S, \phi}(\epsilon) < B_{S,\phi}|\ln(\epsilon)| \epsilon^{1 \over d} \eqno (1.12)$$
One has some analogous lower bounds for $I_{|S|,\phi}(\epsilon)$. They are sharp in cases 1 and 2, and
almost sharp in case 3. (Sharp lower bounds do not hold in general in case 3, as explicit examples show).

\noindent {\bf Theorem 1.3b.} Suppose we are in the setting of Theorem 1.3a). Suppose also that $\phi(0,0)
\neq 0$.

\noindent In case 1 there exists a $A_{S,\phi} > 0$ such that for $\epsilon$ sufficiently small we have
$$I_{|S|, \phi}(\epsilon) > A_{S,\phi}\epsilon^{1 \over d} \eqno (1.13a)$$
In case 2, one similarly has 
$$I_{|S|, \phi}(\epsilon) > A_{S,\phi}|\ln(\epsilon)| \epsilon^{1 \over d} \eqno (1.13b)$$
In case 3, one has analogous almost-sharp lower bounds, at least if $\phi(x,y)$ is nonnegative.
Namely, for any $\delta > 0$ one has
$$I_{|S|, \phi}(\epsilon) > A_{S,\phi, \delta} \epsilon^{{1 \over d} + \delta} \eqno (1.13c)$$

The next lemma will be used to show that the three cases of superadapted coordinates are 
mutually exclusive.

\noindent {\bf Lemma 1.4.} Assume $S(x,y)$ is smooth and is in case 3 of superadapted coordinates 
with $d > 1$. Then for any $M$ one can find a smooth function $S_M(x,y)$
such that $S_M - S$ has a zero of order at least $M$ at the origin, but such that
in a small enough neighborhood $U$ of the origin one has
$$\int_U |S_M(x,y)|^{-{1 \over d}}\,dx \,dy < \infty \eqno (1.14)$$

\noindent {\bf Proof.} Suppose for example that the bisectrix intersects $N(S)$ in the interior of 
the horizontal ray. Then $S_M(x,y) = S(x,y) + x^M$ agrees with $S(x,y)$ to order $M$ at the origin and has 
Newton distance less than that of $S(x,y)$. The Newton distance is also
greater than 1 for $M$ large and the relevant polynomials $(S_M)_e(1,y)$ and $(S_M)_e(-1,y)$ have zeroes 
of order at most 1. Hence $S_M(x,y)$ is in case 1 superadapted coordinates and one 
can apply Theorem 1.1a to conclude that $(1.14)$ is finite for a small enough neighborhood $U$ of the origin. 

\noindent {\bf Lemma 1.5.} Any smooth degenerate phase $S(x,y)$ can be put in superadapted 
coordinates in exactly one of cases 1, 2, or 3.

\noindent {\bf Proof:} In section 7 we will show that one can always put $S(x,y)$ into some superadapted
coodinate system. Equations $(1.11)$ (for $S(x,y)$ and $-S(x,y)$) and $(1.13b)$
cannot simultaneously hold, so a case 2 coordinate system cannot be put in a case 1 or 3 
coordinate system. Suppose $S(x,y)$ has a case 3 
coordinate system as well as a case 1 coordinate system; we will derive a contradiction. Since
it has a case 3 coordinate system, we may adjust $S(x,y)$ to arbitrarily high order and cause $(1.14)$ 
to hold. In its case 1 coordinates, an adjustment of high enough order will not affect the 
fact that it is in case 1 and thus equation $(1.11)$ will still hold. Using the 
well-known relationship between $L^p$ norms and distribution functions (on the function
$(I_{|S|, \phi}(\epsilon))^{-1}$), one then gets that the integral $(1.14)$ in the new coordinates is 
infinite, a contradiction. Thus the three cases are mutually exclusive.

\noindent For oscillatory integrals with smooth phase, one has some analogues of Theorem 1.3. It should be
noted that Theorem 1.6a can be proved using the results of [IKM].

\noindent {\bf Theorem 1.6a.} Suppose $S(x,y)$ is smooth and is in superadapted coordinates with $d > 1$
In cases 1 and 3 as $\lambda \rightarrow \infty$ one has
$$ |J_{S,\phi}(\lambda)| < C\lambda^{-{1 \over d}} \eqno (1.15a)$$
In case 2 one has
$$ |J_{S,\phi}(\lambda)| < C\lambda^{-{1 \over d}}\ln(\lambda)\eqno (1.15b)$$
\noindent {\bf Theorem 1.6b.} Suppose $\phi(x,y)$ is nonnegative with $\phi(0,0) > 0$. 

\noindent In case 1, one has 
$$\limsup_{\lambda \rightarrow \infty} \Big\vert{J_{S,\phi}(\lambda) \over 
\lambda^{-{1 \over d}}}\Big\vert > 0 \eqno (1.16a)$$
In case 2, one has
$$\limsup_{\lambda \rightarrow \infty} \Big\vert{J_{S,\phi}(\lambda) \over 
\lambda^{-{1 \over d}}\ln(\lambda)}\Big\vert > 0 \eqno (1.16b)$$
In case 3, for any $\delta > 0$ one has
$$\limsup_{\lambda \rightarrow \infty} \Big\vert{J_{S,\phi}(\lambda) \over 
\lambda^{-{1 \over d} - \delta}}\Big\vert = \infty \eqno (1.16c)$$
Although we will not prove it here, it can be shown with some additional argument that the conditions
on $\phi(x,y)$ in $(1.16a)-(1.16b)$ can be weakened to just that $\phi(0,0) \neq 0$.

In all of cases 1, 2, and 3 we will divide the domain of integration of the expressions $(1.1)$ and $(1.3)$ 
for $J_{S,\phi}$ and $I_{S,\phi}$ into 4 parts, depending on whether or not $x$ and $y$ are positive or
negative. Adding the resulting formulas and estimates will give the theorems. Without loss of generality
we will always focus on the $x, y > 0$ as the other quadrants are always dealt with the same way. Hence
our goal is to understand $I_{S,\phi}^+$ and $J_{S,\phi}^+$, where
$$I_{S, \phi}^+(\epsilon) =  \int_{\{(x,y): x > 0, y > 0,\,\,0 < S(x,y) < \epsilon\}} \phi(x,y) 
\,dx\,dy \eqno (1.17a)$$
$$J_{S, \phi}^+(\lambda) = \int_{\{(x,y): x > 0, y > 0\}} e^{i \lambda S(x,y)} \phi(x,y)\,dx\,dy\eqno (1.17b)$$
In turn, the domains of $(1.17a)-(1.17b)$ will be written as the union of various ``curved triangles''
(such as those of Lemma 2.0 below). On a given curved triangle, one typically Taylor expands $S(x,y)$ 
or one of its derivatives and then uses Van der Corput-type lemmas in the $x$ or $y$ direction to get 
a desired estimate. For the oscillatory integrals, the traditional van der Corput (see Ch 8 of [S]) is 
used, while for sublevel integrals the version of [C] is used. Van der Corput-type lemmas 
have been considered in some detail, such as in [ArKaCu] and [CaCW], as well as the early work of Vinogradov [Vi]. We refer
to [CaCW] for further results and references. 

Throughout this paper, we will often have a constant $C$ appearing on the right-hand
side of an inequality. This always denotes a constant depending on $S$ and $\phi$. Occasionally
we will need further constants $C'$, $C''$, etc which also depend on $S$ and $\phi$.

\noindent {\bf 2. Some useful lemmas for Cases 1 and 2.} 

Suppose $G$ is an open subset of $\R^2$. Then throughout the course of this paper we will make frequent
use of $I_{S,\phi}^G(\epsilon)$ and $J_{S,\phi}^G(\lambda)$ defined by 
$$I_{S,\phi}^G(\epsilon) = \int_{\{(x,y) \in G: 0 < S(x,y) < \epsilon\}} \phi(x,y) \,dx\,dy
\eqno (2.1a)$$
$$J_{S,\phi}^G(\lambda) = \int_G e^{i \lambda S(x,y)} \phi(x,y)\,dx\,dy \eqno (2.1b)$$
A certain type of $G$ comes up in several contexts in this paper, and relevant estimates we need for 
$I_{S,\phi}^G(\epsilon)$ and $J_{S,\phi}^G(\lambda)$ are given by the following lemma.

\noindent {\bf Lemma 2.0.} Suppose for some $A, m > 0$, $0 < \delta < 1$, we let $G = \{(x,y) 
\in [0,\delta] \times [0,\delta]: 0 <
y < Ax^m\}$, and suppose there are nonnegative integers $a$ and $b$ with $a > b$ and $a \geq 2$ such
that for some constant $C_0$ the following holds on $G$.
$$\partial_y^bS(x,y) > C_0 x^a \eqno (2.2a)$$
If $b = 1$, assume also that
$$\partial_y^2S(x,y) < C_0 x^{a - m} \eqno (2.2b)$$
If $b = 0$, instead of $(2.2b)$ assume also that $a > m + 1$ and that for some constant $C_1$ we have
$$ \partial_xS(x,y) > C_1 x^{a - 1},\,\,\,\,\,\, \partial_x^2S(x,y) < C_0 x^{a - 2} \eqno (2.2c)$$
Then for some $\zeta_{ab} > 0$ and $C$ depending on $S$, $\phi$, and $C_0$ (and $C_1$ if $b = 0$), if the support of 
$\phi$ is contained in $[-\delta,\delta] \times [-\delta,\delta]$ one has 
$$|I_{S,\phi}^G(\epsilon)| < C\epsilon^{m+1 \over a + mb} A^{\zeta_{ab}} \eqno (2.3a)$$
$$|J_{S,\phi}^G(\lambda)| < C|\lambda|^{-{m+1 \over a + mb}} A^{\zeta_{ab}} \eqno (2.3b)$$
\noindent {\bf Proof.} We first consider the case where $b > 0$; we will do the $b = 0$ argument 
afterwards. By $(2.2a)$ and the Van der Corput lemma in 
the $y$ direction (see [C] for example) one has that for a given $x$ we have
$$\big|\{y: |S(x,y)| < \epsilon\}\big| < C \epsilon^{{1 \over b}} x^{-{a \over b}} \eqno (2.4)$$
As a result, if $G^x$ denotes the vertical cross section of $G$ at $x$, of length $Ax^m$, then one has 
$$\big|\{y \in G^x: |S(x,y)| < \epsilon\}\big| < C\min(Ax^m, \epsilon^{{1 \over b}} x^{-{a \over b}})
\eqno (2.5)$$
Consequently, we have
$$|I_{S,\phi}^G| < C \int_0^{\delta} \min(Ax^m, \epsilon^{{1 \over b}} x^{-{a \over b}})\,dx 
\eqno (2.6)$$
It is natural to break the integral $(2.6)$ into two parts, depending on whether or not $Ax^m > 
\epsilon^{{1 \over b}} x^{-{a \over b}}$. The two quantities are equal at $x_0 = A^{-{b \over a + mb}}
\epsilon^{1 \over a + mb}$. The left integral becomes 
$$\int_0^{x_0} Ax^m \,dx = {A \over m+1}x_0^{m+1} = {1 \over m+1} A^{a - b \over a + mb} \epsilon^
{m+1 \over a + mb}\eqno (2.7)$$
The right integral is computed to be
$$\epsilon^{{1 \over b}}\int_{x_0}^{\delta} x^{-{a \over b}}\,dx < \epsilon^{{1 \over b}}\int_{x_0}^{\infty}
x^{-{a \over b}}\,dx $$
$$= {b \over a - b}\epsilon^{1\over b}x_0^{a - b \over b}  = {b \over a - b}A^{a - b \over a + mb} 
\epsilon^{m+1 \over a + mb}\eqno (2.8)$$
Adding together, we obtain that $|I_{S,\phi}^G| < C A^{a - b \over a + mb} \epsilon^{m+1 \over a + mb}$
as needed. 

The estimates for $J_{S,\phi}^G(\lambda)$ for $b \geq 2$ are done in a similar fashion. First 
suppose $b \geq 2$. Then
one can use the usual Van der Corput lemma (see [S] ch 8) in the $y$ direction to obtain
$$|\int_{G^x} e^{i \lambda S(x,y)}\phi(x,y) \,dy| < C|\lambda|^{-{1 \over b}} x^{-{a \over b}}
\eqno (2.9)$$
This is the analogue to $(2.4)$ with $\epsilon$ replaced by $|\lambda|^{-1}$. As a result, similar
to $(2.6)$ we get
$$|J_{S,\phi}^G| < C   \int_0^{\delta} \min(Ax^m, |\lambda|^{-{1 \over b}} x^{-{a \over b}})\,dx  \eqno 
(2.10)$$
The result is  
$$|J_{S,\phi}^G| \leq CA^{a - b \over a + mb} |\lambda|^{-{m+1 \over a + mb}} \eqno (2.11)$$
This gives $(2.3b)$. We next prove $(2.3b)$ when $b = 1$. If one integrates by parts
in the $y$ variable one gets several terms each of which can be bounded using $(2.2a)$ and $(2.2b)$.
If one works it out, one gets that these terms are bounded by $C|\lambda|^{-1}x^{-(a + m)}$. It 
is thus
natural to split the integral into two parts at the $x_0$ satisfying $|\lambda|^{-1}{x_0}^{a + m} = 1$, in 
other words, at $x_0 = |\lambda|^{-{1 \over a + m}}$. We accordingly write $G = G_1 \cup G_2$, with
$G_1$ the portion where $x_0 < |\lambda|^{-{1 \over a + m}}$. We then have
$$J_{S,\phi}^{G_1}(\lambda) < C|G_1| < C'A |\lambda|^{-{m + 1 \over m + a}}\eqno (2.12)$$
As for $J_{S,\phi}^{G_2}(\lambda)$, we integrate by parts in $y$, obtaining the $C|\lambda|^{-1}
x^{-(a + m)}$ factor, and one gets that
$$|J_{S,\phi}^{G_2}(\lambda)| < C\int_{G_2} |\lambda|^{-1}x^{-(m + a)} $$
$$ \leq C\int_{|\lambda|^{-{1 \over a}}}^1 \int_0^{Ax^m}|\lambda|^{-1} x^{-(m + a)} \,dy\,dx $$
$$ = CA\int_{|\lambda|^{-{1 \over a + m}}}^1 |\lambda|^{-1}x^{-a}\,dx $$
$$< CA|\lambda|^{-{m +1 \over m + a}} \eqno (2.13)$$
Adding $(2.12)$ to $(2.13)$ give the oscillatory integral estimates for $b = 1$.

We now consider the oscillatory integral when $b = 0$. Here we do the integrations by parts in the $x$ 
direction. This time, by $(2.2a)$ and $(2.2c)$ an integration by parts incurs a factor of 
$|\lambda|^{-1} x^{-a}$. Hence we 
subdivide $G = G_1 \cup G_2$, where $G_1 = \{(x,y) \in G: 0 < x <|\lambda|^{-{1 \over a}}\}$. Note 
that the measure of $G_1$ is $A|\lambda|^{-{m + 1 \over a}}$, so that
$$|J_{S,\phi}^{G_1}(\lambda)| < C A|\lambda|^{-{m + 1 \over a}} \eqno (2.14)$$
For the $G_2$ piece one obtains
$$|J_{S,\phi}^{G_2}(\lambda)| < C \int_G |\lambda|^{-1} x^{-a} = 
\int_{|\lambda|^{-{1 \over a}}}^1 \int_0^{Ax^m}|\lambda|^{-1} x^{-a} \,dy\,dx $$
$$= C |\lambda|^{-1} \int_{|\lambda|^{-{1 \over a}}}^1A x^{m - a}\,dx $$
$$ = C A|\lambda|^{-{m + 1 \over a}} \eqno (2.15)$$
(For the last equality we use the hypothesis that $a > m + 1$).
Adding $(2.14)$ to $(2.15)$ gives the estimate we seek, $(2.3b)$. Lastly, we prove the bounds for
$I_{S,\phi}^G(\epsilon)$ when $b = 0$. In this case, since $|S(x,y)| > C_0x^a$, we have 
$$|I_{S,\phi}^G| < C|\{(x,y) \in G: C_0x^a < \epsilon\}| < C'\{(x,y) \in G: x < \epsilon^{1 \over a}\}
= C''A\epsilon^{m + 1 \over a} \eqno (2.16)$$
This concludes the proof of Lemma 2.0.

In Lemmas 2.1 and 2.2 below, $S(x,y)$ is a smooth phase function in Case 1 or 2 of superadapted 
coordinates with $d > 1$. If
there is a compact edge $E$ of $N(S)$ such that the bisectrix contains either the upper vertex of 
$E$ or an interior point of $E$, then we denote the equation of this edge by $t_1 + mt_2 = \alpha$,
and for some large but fixed number $N$ we let $A_1 = \{(x,y) \in [0,1] \times [0,1]: y <
{1 \over N} x^m\}$. Similarly, if there is some edge $E'$ with equation $t_1 + m't_2 = \alpha'$ such
that the bisectrix contains either the lower vertex of $E'$ or an interior point of $E'$, we let
$A_2 = \{(x,y) \in [0,1] \times [0,1]: x < {1 \over N} y^{1 \over m'}\}$. Note that in case 1 both $A_1$ and $A_2$ exist. 
We focus our attention on $I_{S,\phi}^{A_i}(\epsilon)$ and $J_{S,\phi}^{A_i}(\lambda)$. The relevant 
information about them (if they exist) is provided by the following lemma.

\noindent {\bf Lemma 2.1.} There exists an $\eta > 0$ such that if the support of $\phi(x,y)$ is
sufficiently small, then for $i = 1, 2$ we have
$$|I_{S,\phi}^{A_i}(\epsilon)| < C\epsilon^{1 \over d} N^{- \eta} \eqno (2.17)$$
$$|J_{S,\phi}^{A_i}(\lambda)| < C|\lambda|^{-{1 \over d}} N^{- \eta} \eqno (2.18)$$
\noindent {\bf Proof.} By symmetry, it suffices to prove the bounds for $A_1$. Let $(a,b)$ denote
the lowest vertex of $E$. Thus $b < a$. Since $(a,b)$ and $(d,d)$ are both on $E$, we have $a + mb = 
(1 + m)d$ or ${1 \over d} = {m + 1 \over a + mb}$. We will show that
the hypotheses of Lemma 2.0 hold for these values of $a$, $b$, and $m$, setting $A = N^{-1}$. Since 
${1 \over d} = {m + 1 \over a + mb}$, Lemma 2.1 will follow.
For a large but fixed $M$, we write $S(x,y) = \sum_{p < M,\,\,\,q < M}s_{pq}x^py^q + E_M(x,y)$. By 
standard estimates, for $0 \leq \alpha, \beta \leq M$ we have
$$|\partial_x^\alpha \partial_y^\beta E_M(x,y)| < C(|x|^{M - \alpha} + |y|^{M - \beta}) \eqno (2.19)$$
  We can write 
$$\partial_y^b S(x,y) = \sum_{p < M,\,\,\,q < M - b}s_{pq}'x^py^q + \partial_y^b E_M(x,y) \eqno (2.20)$$
Here $s_{a0}' \neq 0$. We next show that the sum in $(2.20)$ is dominated by the term $s_{a0}'x^a$ in
a sufficiently small neighborhood of the origin.
To this end, note that $(a,0)$ is a vertex of the Newton polygon of $\partial_y^b S(x,y)$, and that a 
horizontal ray and an edge of this Newton polygon with equation $t_1 + mt_2 = a$ intersect at $(a,0)$. As
a result, for any term $s_{pq}'x^py^q$ in the sum of $(2.20)$ other than $s_{a0}'x^a$, either 
$p \geq a$, or $p < a$, $q > 0$, and $p + mq \geq a$. Correspondingly, we let $T_1 = \{(p,q): 
p < a,\,\,0 < q < M - b,\,\,p + mq \geq a\}$ and $T_2 = \{(p,q): a \leq p \leq M,\,\,0 \leq q 
\leq b - a,\,\,(p,q) \neq (a,0)\}$ and we rewrite $(2.20)$ as
$$\partial_y^b S(x,y) = s_{a0}'x^a + \sum_{T_1}s_{pq}'x^py^q + \sum_{T_2}s_{pq}'x^py^q + 
\partial_y^b E_M(x,y) \eqno (2.21)$$
We examine a given term $s_{pq}'x^py^q$ in the $T_1$ sum. Since $(x,y)$ is in the domain $A_1$, one
has $y < {x^m \over N}$. As a result, $|s_{pq}'x^py^q| < |s_{pq}'N^{-q}x^{p + mq}| \leq |s_{pq}'N^{-q}
x^a|$. Thus so long as $N$ is chosen sufficiently large (depending on $M$), we may assume that the 
absolute value of the $T_1$ sum is at most ${1 \over 4}|s_{a0}'| x^a$. We next examine a term 
$s_{pq}'x^py^q$ in the $T_2$ sum. Here we have $|s_{pq}'x^py^q|< |s_{pq}'x^a|(|x| + |y|)$. Hence if
we are in a sufficiently small neighborhood of the origin, we can assume the absolute value of the 
$T_2$ sum is also at most ${1 \over 4}|s_{a0}'| x^a$. Lastly, by $(2.19)$, in a sufficiently small
neighborhood of the origin $|\partial_y^b E_M(x,y)|$ is also bounded by ${1 \over 4}|s_{a0}'| x^a$.
Putting these together, we conclude that on the domain $A_1$ we have
$$|\partial_y^b S(x,y)| > {1 \over 4}|s_{a0}'|x^a \eqno (2.22)$$
Note that when $b > 1$, $(2.22)$ ensures hypotheses of Lemma 2.0 hold. Thus one may apply Lemma 2.0, 
giving Lemma 2.1. When $b = 0$, in order to apply Lemma 2.0 one needs also 
that $a > m + 1$ and that $(2.2c)$ 
holds. But taking an $x$ derivative of $S(x,y)$ just shifts the Newton polygon to the left by 1, so
exactly as in $(2.22)$ the first and second derivative conditions of $(2.2c)$ will hold. As for the
requirement that $a > m + 1$, note that $(d,d)$ and $(a,0)$ are on the
line $t_1 + mt_2 = \alpha$ and that $d > 1$. This means $a = a + m0 = d + md > 1 + m$ as needed.

Lastly, we show that the supplemental hypothesis $(2.2b)$ holds when $b = 1$. It helps to view things 
in $(x,y')$ coordinates where $(x,y) = (x,x^my')$. The condition $(2.2b)$ becomes that 
$|\partial_{y'}^2S(x,x^my')| \leq  Cx^{a + m}$. 
Also, since the terms $s_{pq}x^py^q$ of $S(x,y)$'s Taylor expansion with minimal $p + mq$ ($= \alpha$)
are exactly the terms of $S_E(x,y)$, the expansion $S(x,y) = \sum_{p < M,\,\,\,q < M}s_{pq}x^py^q +
E_M(x,y)$ becomes of the following form, where $T_M$ is a polynomial in $y'$ and a fractional power of 
$x$.
$$S(x,x^my') = x^{\alpha}S_E(1,y') + x^{\alpha + \epsilon}T_M(x,y') + E_M(x,x^my') \eqno (2.23)$$
Using $(2.23)$ and the error estimates $(2.19)$ we have $|\partial_{y'}^2S(x,x^my')| \leq  Cx^{\alpha}$. 
But $(a,b) = (a,1)$ is on the line $t_1 + mt_2 = 
\alpha$ and therefore  $a + m = \alpha$. This gives the second derivative bounds of $(2.2b)$, and thus
the hypotheses of Lemma 2.0 hold. This concludes the proof of Lemma 2.1.

\noindent {\bf Lemma 2.2.} Let $S(x,y)$ be a smooth phase function in case 1 or 2 superadapted 
coordinates with $d > 1$. Suppose $e$ is a compact edge of $N(S)$ intersecting the bisectrix and
has equation given by $t_1 + mt_2 = 
\alpha$. Define $B = \{(x,y) \in [0,1] \times [0,1]: {x^m \over N} < y < Nx^m\}$, where $N$ is
some large but fixed constant. Then if the support of $\phi(x,y)$ is sufficiently small, depending on 
$N$, one has the estimates
$|I_{S,\phi}^B(\epsilon)| < C_N \epsilon^{1 \over d}$ and $|J_{S,\phi}^B
(\lambda)| < C_N|\lambda|^{-{1 \over d}}$.

\noindent {\bf Proof.} In the above $(x,y')$ coordinates one has
$$I_{S,\phi}^B(\epsilon) = \int_{\{(x,y') \in [0,1] \times [N^{-1},N]: 0 < S(x,x^my') < \epsilon\}}
x^m \phi(x,x^my') \,dx\,dy'$$
$$J_{S,\phi}^B(\lambda) = \int_{[0,1] \times [N^{-1},N]} e^{i \lambda S(x,x^my')} x^m \phi(x,x^my')
\,dx\,dy'$$
In view of $(2.23)$, the zeroes of $S_e(1,y')$ might be expected to play a significant role in the
analysis. To this end, we assume that $N$ is large enough so that any zeroes
of $S_e(1,y')$ for $y' > 0$ are in $[N^{-1},N]$, and denote these zeroes by $z_1,...,z_k$ (if there are
any).  Let $v_1,...v_k$ denote the orders of these zeroes, and let $I_i$ denote the interval 
$[z_i - {1 \over N}, z_i + {1 \over N}]$. By $(2.23)$ and the error
term derivative estimates $(2.19)$ we can assume that on the (sufficiently small, depending on $N$) 
support of $\phi(x,x^my')$, if $(x,y') \in [0,1] \times I_i$ then 
$$|\partial_{y'}^{v_i} S(x,x^my')| \geq  Cx^{\alpha} \eqno (2.24a)$$
In the case that $v_i = 1$, we similarly have
$$|\partial_{y'}^2 S(x,x^my')| \leq  Cx^{\alpha} \eqno (2.24b)$$
We now translate this into the orginal $(x,y)$ coordinates. Each $[0,1] \times I_i$ becomes a set 
$D_i$ of the form
$\{(x,y) \in [0,1] \times [0,1]: (z_i - {1 \over N})x^m < y < (z_i + {1 \over N})x^m$, and on $D_i$ 
$(2.24a)$ becomes
$$|\partial_y^{v_i} S(x,y)| \geq  Cx^{\alpha - mv_i} \eqno (2.25a)$$
In the case of $v_i = 1$, $(2.24b)$ becomes
$$|\partial_y^2 S(x,y)| \leq  Cx^{\alpha - 2m} \eqno (2.25b)$$
Since we are superadapted coordinates, $0 < v_i < d$. Thus $(m + 1)v_i < (m + 1)d$. Since $(d,d)$ is 
on the edge $e$, we have $\alpha = (m + 1)d$. Thus $(m + 1)v_i < \alpha$ or $\alpha - mv_i > v_i$. 
As a result, the sets $D_i$ have 
vertical cross sections of length ${2 \over N} x^m$ and satisfying $(2.25a)-(2.25b)$ with $\alpha - 
mv_i > v_i$. Hence after doing a coordinate change of the form $(x,y) \rightarrow (x,y - f(x))$, we 
are in the set-up of Lemma 2.0 and for some $\eta > 0$ we get 
$$|I_{S,\phi}^{D_i}|  \leq  C\epsilon^{m + 1 \over \alpha} N^{-\eta}\eqno (2.26a)$$
$$|J_{S,\phi}^{D_i}|  \leq C|\lambda|^{-{m + 1 \over \alpha}} N^{-\eta}\eqno (2.26b)$$
Since ${m + 1 \over \alpha} = {1 \over d}$ the above becomes
$$|I_{S,\phi}^{D_i}|  \leq  C\epsilon^{{1 \over d}} N^{-\eta} \eqno (2.27a)$$
$$|J_{S,\phi}^{D_i}|  \leq C|\lambda|^{-{1 \over d}} N^{-\eta}\eqno (2.27b)$$
Next, write $[N^{-1},N] - \cup_i I_i$ as the union of intervals $J_i$. Then since $S_e(1,y')$ has no 
zeroes on any $J_i$, by the expansion $(2.23)$ and the error derivative bounds $(2.19)$, if $\delta_N$
is sufficiently small then on $[0,\delta_N] \times J_i$ we have
$$|\partial_x S(x,x^my')| > C_Nx^{\alpha - 1} \hskip 25 pt |\partial_x^2 S(x,x^my')| < C_N'x^{\alpha - 2}
\eqno (2.28)$$
Separating at $x = |\lambda|^{-{1 \over \alpha}}$ and integrating the right portion by parts in $x$ 
using $(2.28)$ as in the proof of Lemma 2.0 gives
$$|\int_{[0,\delta_N] \times J_i} e^{i \lambda S(x,x^my')} x^m \phi(x,x^my')| \,dx\,dy' < 
C_N''|\lambda|^{-{m + 1 \over \alpha}} \eqno (2.29a)$$
Converting back into $(x,y)$ coordinates and using that ${1 \over d} = {m + 1 \over \alpha}$,
$(2.29a)$ becomes the following, where $E_i$ denotes the set $[0,\delta_N] \times J_i$ in the 
$(x,y)$ coordinates.
$$|\int_{E_i} e^{i \lambda S(x,y)} \phi(x,y)\,dx\,dy| \leq C_N''|\lambda|^{-{1 \over d}} 
\eqno (2.29b)$$
Moving now to $I_{S,\phi}(\epsilon)$, by  $(2.23)$ one has $|S(x,x^my')| > C_Nx^{\alpha}$ on $[0,
\delta_N] \times J_i$ for sufficiently small $\delta_N > 0$ and some $C_N$ (not necessarily the same
constant as above). As a result, 
$$|\int_{\{(x,y') \in [0,\delta_N] \times J_i: 0 < S(x,x^my') < \epsilon\}} x^m \phi(x,x^my') \,dx\,dy'| $$
$$\leq |\int_{\{(x,y') \in [0,\delta_N] \times J_i: 0 < C_Nx^{\alpha} < \epsilon\}} x^m \phi(x,x^my') \,dx\,dy'| $$
$$ < C_N'\int_0^{\epsilon^{1 \over \alpha}} x^m \,dx$$
$$ = C_N'\epsilon^{{m + 1 \over \alpha}} \eqno (2.30a)$$
Again going back to the $(x,y)$ coordinates and using that ${1 \over d} = {m + 1 \over \alpha}$, we 
conclude that
$$|\int_{\{(x,y) \in E_i: 0 < S(x,y) < \epsilon\}} \phi(x,y) \,dx\,dy| \leq C_N\epsilon^{1 \over d} 
\eqno (2.30b)$$
Lemma 2.2 now follows by adding $(2.29b)-(2.30b)$ to $(2.27a)-(2.27b)$.

\noindent {\bf 3. Case 1 proofs.} 

Assume now that we are in Case 1. We start by proving the upper bounds for smooth phase.

\noindent {\bf Theorem 3.1.} The right-hand sides of $(1.11)$ and $(1.15a)$ hold.

\noindent {\bf Proof.} In case 1 of superadapted coordinates, the domain of integration of 
$I_{S,\phi}^+(\epsilon)$ or $J_{S,\phi}^+(\lambda)$ is the union of $A_1$, $A_2$, and $B$, where $A_1$
and $A_2$ are as in Lemma 2.1 and
where $B$ is as in Lemma 2.2 for the edge of $N(S)$ intersecting the bisectrix. Thus
the theorem follows by fixing some $N$ and adding the inequalities of Lemmas 2.1 and 2.2 to the 
corresponding inequalities for the other quadrants.

\noindent The next lemma will be useful in getting the formulas for the real-analytic case.

\noindent {\bf Lemma 3.2.} Suppose $S(x,y)$ is a smooth case 1 phase function like before. Then there
is a natural number $D < d$ and a neighborhood $U$ of the origin such that if $\phi(x,y)$ is supported 
in $U$ and is zero on a neighborhood of the origin, then $|I_{S,\phi}(\epsilon)| < C\epsilon^{1 \over 
D}$. 

\noindent {\bf Proof:} It suffices to fix some $N$ and show that each $I_{S,\phi}^{A_i}(\epsilon)$, 
$I_{S,\phi}^{D_i}(\epsilon)$, $I_{S,\phi}^{E_i}(\epsilon)$ satisfies the upper bounds, where the $D_i$ and $E_i$ correspond
to the edge of $N(S)$ intersecting the bisectrix. We start with $I_{S,\phi}^
{A_i}$. Without loss of generality we may take $i = 1$. Since $\phi(x,y)$ is zero on a neighborhood of the origin, there is some $\delta > 0$ such
that $\phi(x,y)$ is zero for $(x,y) \in A_1$ with $0 < x < \delta$. $(2.22)$ says that 
$\partial_y^eS(x,y)$
is bounded below on $A_1$, where $e < d$ denotes the $y$-coordinate of the lower vertex of the edge
of $N(S)$ intersecting the bisectrix. If $e > 0$, by the Van der Corput lemma in the $y$ direction, we have
$$|\{y \in A_1^x:  0 < S(x,y) < \epsilon\}| < C\epsilon^{{1 \over e}} \eqno (3.1)$$
Thus we have
$$|I_{S,\phi}^{A_1}(\epsilon)| \leq C \int_{\delta}^1 |\{y \in A_1^x:  0 < S(x,y) < \epsilon\}|\,dx < C\epsilon
^{{1 \over e}} \eqno (3.2)$$
This is the desired estimate for $e > 0$. If $e = 0$, then $(2.22)$ says that $S(x,y)$ is bounded below 
on the support of the integrand of $I_{S,\phi}^{A_1}(\epsilon)$ and then $(3.2)$ holds trivially. Thus we have 
the desired bounds for the $I_{S,\phi}^{A_1}(\epsilon)$.
 The $I_{S,\phi}^{D_i}(\epsilon)$ are dealt with in a similar way. This time, one uses $(2.24a)$  
to obtain $|I_{S,\phi}^{D_i}(\epsilon)| < C\epsilon^{{1 \over v_i}}$.

Lastly, we look at the $|I_{S,\phi}^{E_i}(\epsilon)|$.
As mentioned below $(2.29b)$, by $(2.23)$ one has $|S(x,y)| > Cx^{\alpha}$ on each $E_i$.
Hence since $\phi(x,y)$ is zero on a neighborhood of the origin, $|S(x,y)|$ is bounded below on the
support of the integrand of $I_{S,\phi}^{E_i}(\epsilon)$ and $(3.2)$ again holds trivially. This completes
the proof of Lemma 3.2.

We now move to the case 1 formulas of Theorems 1.1. We write $\phi =
\phi_1 + \phi_2$, where $\phi_1$ is supported in a smaller neighborhood of the origin and $\phi_2$ is 
zero on a neighborhood of the origin. By Lemma 2.3, $\lim_{\epsilon \rightarrow 0} {I_{S,\phi_2}
(\epsilon) \over \epsilon^{1 \over d}} = 0$ and therefore $\lim_{\epsilon \rightarrow 0} {I_{S,\phi}
(\epsilon) \over \epsilon^{1 \over d}} =$  $\lim_{\epsilon \rightarrow 0} {I_{S,\phi_1}
(\epsilon) \over \epsilon^{1 \over d}}$. Hence when proving the limit $(1.6)$ one can always 
replace $\phi$ by $\phi_1$ at will, regardless of how small the support of $\phi_1$ is.  

Our strategy will involve fixing $N$ and analyzing the $I_{S,\phi_1}^{E_i}$, where $\phi_1$ 
has small support depending on $N$. Lemma 2.1 and $(2.27a)$ will ensure that the contributions
of the $I_{S,\phi_1}^{D_i}$ and $I_{S,\phi_1}^{A_i}$ will be $O(N^{-\eta})$ smaller than that of the
$I_{S,\phi_1}^{E_i}$ as $\epsilon$ goes to zero. 
Letting $N$ go to infinity will give $\lim_{\epsilon \rightarrow 0} { I_{S,\phi_1}^+(\epsilon) \over 
\epsilon^{1 \over d}} = \sum_i \lim_{\epsilon \rightarrow 0} { I_{S,\phi_1}^{E_i}(\epsilon) \over 
\epsilon^{1 \over d}}$. Adding up the latter limits along with their analogues in the 
other three quadrants will give  $(1.6)$. 

\noindent {\bf Proof of $(1.6)$.} It suffices to assume $\phi(0,0) > 0$ as the general case can be
obtained by writing $\phi = \phi' - \phi''$ where $\phi'(0,0)$ and  $\phi''(0,0)$ are positive.
We work in the $(x,y')$ coordinates like above. Namely, we write
$$I_{S,\phi}^{E_i}(\epsilon) = \int_{\{(x,y') \in [0,1] \times J_i: 0 < S(x,x^my') < \epsilon\}} x^m \phi(x,x^my')
\,dx\,dy' \eqno (3.3)$$
$J_i$ was defined so that $S_0(1,y')$ has no zeroes on $J_i$. As a result, by $(2.23)$ (and using 
$(2.19)$ to deal with $E_M(x,x^my')$), we 
may let $\delta > 0$ such that on $[0,\delta]\times J_i$ either $S(x,x^my')$ is negative, or we have
$$0 < (1 - {1 \over N})S(x,x^my') \leq x^{\alpha}S_0(1,y') \leq (1 + {1 \over N})S(x,x^my') \eqno (3.4)$$
Shrinking $\delta$ further if necessary, we assume $\delta$ is small enough so that on $[0,\delta]
\times J_i$ we have
$$ \phi(x,x^my') < (1 + {1 \over N})\phi(0,0) \eqno (3.5)$$
As described above, one can multiply $\phi(x,y)$ by a cutoff function supported on $|x| < \delta$
without affecting $\lim_{\epsilon \rightarrow 0} { I_{S,\phi_1}^{E_i}(\epsilon) \over 
\epsilon^{1 \over d}}$. Hence we assume $\phi(x,y)$ is supported on $|x| < \delta$, and that
the multiplying cutoff was chosen so that for $(x,y') \in [0,\delta] \times J_i$, equation
$(3.5)$ holds. We also assume the multiplying cutoff was chosen so that for $(x,y') \in 
[0,{\delta \over 2}] \times J_i$ one has 
$$(1 - {1 \over N})\phi(0,0) < \phi(x,x^my') \eqno (3.6)$$
We now proceed to our main estimates. If $S(x,x^my')$ is negative on $[0,\delta]\times J_i$, 
$I_{S,\phi}^{E_i}(\epsilon)$ becomes zero. If on the other hand $(3.4)$ holds, then by $(3.4)$ and 
$(3.5)$ one has
$$I_{S,\phi}^{E_i}(\epsilon) \leq (1 + {1 \over N})\phi(0,0) \int_{\{(x,y') \in [0,\delta] \times J_i: 0 < x^{\alpha}
S_0(1,y') <(1 + {1 \over N}) \epsilon\}} x^m  \,dx\,dy' \eqno (3.7a)$$
On the other hand, by $(3.4)$ and $(3.6)$ we also have
$$I_{S,\phi}^{E_i}(\epsilon) \geq (1 - {1 \over N})\phi(0,0) \int_{\{(x,y') \in [0,{1 \over 2}\delta] \times J_i: 
0 < x^{\alpha} S_0(1,y') < (1 - {1 \over N}) \epsilon\}} x^m  \,dx\,dy' \eqno (3.7b)$$
We now change coordinates from $x$ to $x' = x^{m+1}$ in the integrals $(3.7a)-(3.7b)$. Using the fact
that ${\alpha \over m + 1} = d$, which follows from the fact that $(d,d)$ is on the line $t_1 + mt_2 = 
\alpha$, $(3.7a)-(3.7b)$ become the following, where $\delta' = \delta^{m+1}$.
$$I_{S,\phi}^{E_i}(\epsilon) \leq (m+1)^{-1}(1 + {1 \over N})\phi(0,0) \int_{\{(x,y') \in [0,\delta '] \times J_i: 0 <
x^d S_0(1,y') <(1 + {1 \over N}) \epsilon\}}  \,dx\,dy' \eqno (3.8a)$$
$$I_{S,\phi}^{E_i}(\epsilon) \geq (m+1)^{-1}(1 - {1 \over N})\phi(0,0) \int_{\{(x,y') \in [0,{1 \over 2}\delta '] \times J_i: 
0 < x^d S_0(1,y') < (1 - {1 \over N}) \epsilon\}} \,dx\,dy' \eqno (3.8b)$$
Doing the $x$ integrals first, $(3.8a)-(3.8b)$ become
$$I_{S,\phi}^{E_i}(\epsilon) \leq (m+1)^{-1}(1 + {1 \over N})\phi(0,0) \int_{J_i} \min (\delta ', 
[(1 + {1 \over N})\epsilon]^{1 \over d} S_0(1,y')^{-{1 \over d}})\,dy' \eqno (3.9a)$$
$$I_{S,\phi}^{E_i}(\epsilon) \geq (m+1)^{-1}(1 - {1 \over N})\phi(0,0) \int_{J_i} \min ({\delta '
\over 2}, [(1 - {1 \over N}) \epsilon]^{1 \over d} S_0(1,y')^{-{1 \over d}})\,dy'\eqno (3.9b)$$
These can be written as 
$${I_{S,\phi}^{E_i}(\epsilon) \over \epsilon^{1 \over d}} \leq (m+1)^{-1}(1 + {1 \over N})^{d + 1 \over d} \phi(0,0)
\int_{J_i}\min ({\delta ' \over [(1 + {1 \over N})\epsilon]^{1 \over d}}, S_0(1,y')^{-{1 \over d}})\,dy'
\eqno (3.10a)$$
$${I_{S,\phi}^{E_i}(\epsilon) \over \epsilon^{1 \over d}} \geq (m+1)^{-1}(1 - {1 \over N})^{d + 1 \over d} \phi(0,0)
\int_{J_i}\min ({\delta ' \over 2[(1 - {1 \over N})\epsilon]^{1 \over d}}, S_0(1,y')^{-{1 \over d}})\,dy'
\eqno (3.10b)$$
We now take limits of $(3.10a)-(3.10b)$ as $\epsilon \rightarrow 0$. We obtain
$$ \liminf_{\epsilon \rightarrow 0} {I_{S,\phi}^{E_i}(\epsilon) \over \epsilon^{1 \over d}} \geq 
(m+1)^{-1}(1 - {1 \over N})^{d + 1 \over d} \phi(0,0) \int_{J_i} S_0(1,y')^{-{1 \over d}}\,dy' 
\eqno (3.11a)$$
$$\limsup_{\epsilon \rightarrow 0} {I_{S,\phi}^{E_i}(\epsilon) \over \epsilon^{1 \over d}}
\leq (m+1)^{-1}(1 + {1 \over N})^{d + 1 \over d} \phi(0,0) \int_{J_i} S_0(1,y')^{-{1 \over d}}\,dy' \eqno (3.11b)$$
Note that the integrals here are automatically finite since we are in superadapted coordinates and 
therefore all zeroes of $S_0(1,y')$ are of order at most $d - 1$. We now take limits 
as $N \rightarrow \infty$. If the endpoints of $J_i$ are two zeroes $z$ and $z'$
of $S_0(1,y')$, then the interval will converge to $[z,z']$. Otherwise, the left endpoint of $J_i$ may
converge to zero and the right endpoint may go off to $\infty$. In any event, the $J_i$ goes to 
some (possibly unbounded) interval $K_i$, and $(3.11a)-(3.11b)$ both converge to 
$$ (m+1)^{-1} \phi(0,0) \int_{K_i} S_0(1,y')^{-{1 \over d}}\,dy' \eqno (3.12)$$
The integral in $(3.12)$ is finite since $S_0(1,y')$ must have degree greater than $d$ and has no zeroes
of order $d$ or greater since the coordinate system being used is superadapted.

Furthermore, as $N \rightarrow \infty$, the upper bounds of Lemma 2.1 and $(2.27a)$ for $I_{S,\phi}^{D_i}(\epsilon)$ and 
$I_{S,\phi}^{A_i}(\epsilon)$ go to zero. There is no issue of the constants appearing depending on $N$
due to the cutoffs' dependence on $N$; the bounds of Lemma 2.1 and $(2.27a)$ depend on 
$||\phi||_{\infty}$ and no other properties of $\phi$. Hence we conclude that $\lim_{\epsilon 
\rightarrow 0} { I_{S,\phi_1}^+(\epsilon) \over \epsilon^{1 \over d}}$ is
given by the sum of $(3.12)$ over all $i$, along with their analogues in the other three quadrants.
This gives exactly $(1.6)$ (recall the terms where $S(x,x^my')$ is negative gives no contribution)
and we are done.

Lastly, we prove $(1.13a)$ of Theorem 1.3a. Either $S_0(1,y)$ is positive for some
$y > 0$ or $-S_0(1,y)$ is positive for some $y > 0$.
Since the result we are trying to prove is symmetric in $S$ and $-S$, without loss of generality
we may assume that $S_0(1,y)$ has this property. Then $(1.6)$ says that  $I_{S,\phi}(\epsilon) > 
A_{S,\phi}\epsilon^{1 \over d}$ for sufficiently small $\epsilon$. This implies $(1.13a)$ since 
$I_{|S|,\phi}(\epsilon) \geq I_{S,\phi}(\epsilon)$ and we are done.

\noindent {\bf 4. Case 2 proofs.} 

We now assume $S(x,y)$ is a smooth phase function in case 2 of
superadapted coordinates with $d > 1$. Hence the bisectrix intersects $N(S)$ at $(d,d)$.  
As in the case 1 proofs, we consider $I_{S,\phi}^+(\epsilon)$ and 
$J_{S,\phi}^+(\lambda)$ and we will do some subdivisions of the domains of these integrals to prove
the estimates and formulas. We proceed as follows. If $(d,d)$ is the lower vertex of a
compact edge $e_2$, denote its equation by $t_1 + m_2t_2 = \alpha_2$. Where $N$ is a 
large natural number, fixed for now, we let $A_2 = 
\{(x,y) \in [0,1]  \times [0,1]: x < {1 \over N}y^{1 \over m_2}\}$ and $B_2 = \{(x,y) \in [0,1] \times 
[0,1]: {1 \over N}y^{1 \over m_2} < x < N y^{1 \over m_2}\}$. Define $C_2 = A_2 \cup B_2$.
If $(d,d)$ is not the lower vertex of a compact edge (i.e. $(d,d)$ is on the vertical ray), then 
define $C_2 = \{(x,y) \in [0,1]  \times [0,1]: x < y^L\}$. Here $L$ is large number to be determined 
by our future arguments. Similarly, if there is 
a compact edge $e_1$ whose upper vertex is $(d,d)$, we write its equation as $t_1 + m_1t_2 = \alpha_1$.
We then let $A_1 = \{(x,y) \in [0,1] \times [0,1]: y < {1 \over N}x^{m_1}\}$ and  $B_1 = \{(x,y) \in [0,1] 
\times [0,1]: {1 \over N}x^{m_1} < y < N x^{m_1}\}$.  We then define $C_1 = A_1 \cup B_1$. If $(d,d)$ is
not the upper vertex of a compact edge, then define $C_1 = \{(x,y) \in [0,1]  \times [0,1]: y < x^L\}$.

In all cases, define  $D = [0,1]\times [0,1] - (A_2 \cup B_2)$
We will see that the terms $I_{S,\phi}^D(\epsilon)$ and $J_{S,\phi}^D(\lambda)$ dominate; the contributions
from $C_1$ and $C_2$ to the main term of the asymptotics can be made arbitrarily small as $L \rightarrow
\infty$. In fact,
for the case when $C_i$ comes from a compact edge, Lemmas 2.1 and 2.2 give that 
$$|I_{S,\phi}^{A_i}(\epsilon)| < CN^{-{\eta}}\epsilon^{1 \over d}\hskip 43pt 
|I_{S,\phi}^{B_i}(\epsilon)| < C_N\epsilon^{1 \over d} \eqno (4.1a)$$
$$|J_{S,\phi}^{A_i}(\lambda)| < CN^{-{\eta}}|\lambda|^{-{1 \over d}}\hskip 43pt
 |J_{S,\phi}^{B_i}(\lambda)| < C_N|\lambda|^{-{1 \over d}}\eqno (4.1b)$$
 \noindent {\bf Lemma 4.1.} For sufficiently small $\epsilon$ we have 
$$|I_{S,\phi}^D(\epsilon)| < C|\ln(\epsilon)|\epsilon^{1 \over d} \hskip 44 pt |J_{S,\phi}^D(\lambda)|
< C\ln|\lambda||\lambda|^{-{1 \over d}} \eqno (4.2a)$$
Furthermore, for $i = 1,2$ if $C_i$ derives from the horizontal or vertical ray we have
$$|I_{S,\phi}^{C_i}(\epsilon)| < {C \over L + 1}|\ln(\epsilon)|\epsilon^{1 \over d} \hskip 44 pt 
|J_{S,\phi}^{C_i}(\lambda)| < {C \over L + 1}\ln|\lambda||\lambda|^{-{1 \over d}} \eqno (4.2b)$$ 

\noindent {\bf Proof:} We start with the proof of $(4.2a)$. We divide $D = D_1 \cup D_2$, where 
$D_1 = \{(x,y) \in D: y < x^m\}$. Here $m$ is chosen such that if $e_1$ exists, then
$m < m_1$, and if $e_2$ exists then $m > m_2$. The estimates for $D_1$ and $D_2$ are proven the same
way, so we
restrict our attention to proving the estimates for $D_1$. Taylor expand $S(x,y)$ arount the origin as 
$S(x,y) = \sum_{a < M,\,\,\,b < M}s_{ab}x^ay^b + E_M(x,y)$, where like in $(2.19)$ for
$0 \leq \alpha, \beta \leq M$ we have
$$|\partial_x^\alpha \partial_y^\beta E_M(x,y)| < C(|x|^{M - \alpha} + |y|^{M - \beta}) \eqno (4.3)$$
Correspondingly, we have
$$\partial_y^d S(x,y) = \sum_{a < M,\,\,\,b < M - d}s_{ab}'x^ay^b + {\partial_y^d}E_M(x,y) \eqno 
(4.4)$$
The Newton polygon of $\partial_y^d S(x,y)$ has a vertex at $(d,0)$ and thus $s_{d0}' \neq 0$.
If $(d,d)$ is the lower vertex of some compact edge $e_2$, then there is a compact edge $e'$ of the
Newton polygon of 
$\partial_y^d S(x,y)$ containing $(d,0)$ with equation $t_1 + m_2t_2 = d$. Hence every nonzero
$s_{ab}'x^ay^b$ appearing in $(4.4)$ satisfies $a + m_2b \geq d$ and we can rewrite $(4.4)$ as 
$$\partial_y^d S(x,y) =  s_{d0}'x^d + \sum_{a < d,\,\,\,0 < b < M - d,\,\,\,a + m_2b \geq d}s_{ab}'
x^ay^b$$
$$ + \sum_{d \leq a < M,\,\,\,0 \leq b < M - d,\,\,\,(a,b) \neq (d,0)}s_{ab}'x^ay^b + {\partial_y^d}
E_M(x,y) \eqno (4.5)$$
If $(d,d)$ is not the lower vertex of such a compact edge, $(4.5)$ is still valid if we take the first
sum to be empty. If the first sum is not empty, then since $y < N^{-m_2}x^{m_2}$ for all $(x,y) \in D$,
if $(x,y)$ is in $D$ then each term $s_{ab}'x^ay^b$ in the first
sum of $(4.5)$ is bounded in absolute value by $N^{-m_2b}|s_{ab}'|x^{a + m_2b} \leq 
N^{-m_2b}|s_{ab}'|x^d$. Thus if $N$ were chosen sufficiently large, the absolute value 
of the whole first sum is less than ${1 \over 4}|s_{d0}'|x^d$. 

Next, note that the absolute value of a given term 
$s_{ab}'x^ay^b$ in the second sum is at most $|s_{ab}'|x^d(x + y)$. As a result, if the support of
$\phi(x,y)$ is sufficiently small, then for $(x,y)$ in this support the absolute value of the second
sum is also at most ${1 \over 4}|s_{d0}'|x^d$. Similarly, if the support of $\phi(x,y)$ is sufficiently 
small, then by $(4.2)$ and the fact that $y < x^m$, for $(x,y)$ in this support
$|{\partial_y^d}E_M(x,y)|$ can also be assumed to be at most ${1 \over 4}|s_{d0}'|x^d$. Consequently, for
such $(x,y)$ in the support of $\phi$ we can assume
$$|\partial_y^d S(x,y)| > {1 \over 4}|s_{d0}'|x^d \eqno (4.6)$$
Denote the vertical cross section of $D$ at $x$ by $D_x$. By $(4.6)$ and the measure version of 
Van der Corput's lemma, for each $x$ in this range we have
$$|\{y \in D_x : 0 < S(x,y) < \epsilon\}| < C\epsilon^{1 \over d}{1 \over x}$$
Also, $|\{y \in D_x: 0 < S(x,y) < \epsilon\}|$ is at most $|D_x|$, which is at most 
$x^m$ if the support of $\phi$ is sufficiently small, which we may assume. Thus we have
$$|\{y \in D_x : 0 < S(x,y) < \epsilon\}| <  C\min(x^m,\epsilon^{1 \over d}{1 \over x})$$
Consequently, integrating with respect to $y$ first one has
$$|I_{S,\phi}^D(\epsilon)| < C\int_0^1 \min(x^m,\epsilon^{1 \over d}{1 \over x})
\,dx \eqno (4.7)$$
For small enough $\epsilon$, the quantities $x^m$ and $\epsilon^{1 \over d}{1 \over x}$
are equal at $x_0 = \epsilon^{1 \over (m + 1)d}$, with $x^m$ smaller on the left 
and $\epsilon^{1 \over d}{1 \over x}$ smaller on the right. Doing a computation gives
$$ \int_0^1 \min(x^m,\epsilon^{1 \over d}{1 \over x})\,dx = {1 \over m + 1}\epsilon^
{1 \over d} + {1 \over (m + 1)d}|\ln(\epsilon)|\epsilon^{1 \over d} \eqno (4.8)$$
Hence $|I_{S,\phi}^D(\epsilon)| < C|\ln(\epsilon)|\epsilon^{{1 \over d}}$ as desired. As for the 
$J_{S,\phi}^D(\lambda)$, by $(4.6)$ the traditional Van der Corput lemma in the $y$ direction gives
$$|\int_{D_x} e^{i \lambda S(x,y)} \phi(x,y) \,dy| \leq C |\lambda|^{-{1 \over d}}x^{-1}\eqno (4.9)$$
Consequently, we have
$$|J_{S,\phi}^D(\lambda)| < C\int_0^1 \min(N^{-m}x^m,|\lambda|^{-{1 \over d}}{1 \over x})\,dx
\eqno (4.10)$$
This is exactly $(4.8)$ with $\epsilon$ replaced by $|\lambda|^{-1}$. Thus instead of $(4.8)$ for large
$|\lambda|$ we get the estimate 
$$|J_{S,\phi}^D(\lambda)| < C\ln|\lambda||\lambda|^{-{1 \over d}} \eqno (4.11)$$
This completes the proof of $(4.2a)$. Equation $(4.2b)$ is done the same way; the only difference is that 
$x^m$ is replaced by $x^L$. Equation $(4.8)$ and its oscillatory integral analogue then give $(4.2b)$ 
and we are done.

\noindent We now have proven the upper bounds for the smooth case:

\noindent {\bf Lemma 4.2.} Equations $(1.12)$ and $(1.15b)$ hold.

\noindent {\bf Proof.} Add $(4.2a)$ to $(4.2b)$ or $(4.1a)-(4.1b)$.

\noindent Our next result is an analogue of Lemma 3.2.

\noindent {\bf Lemma 4.3.} There is a neighborhood $U$ of the origin such that if $\phi(x,y)$ is 
supported in $U$ and $\phi(x,y)$ is zero in a neighborhood of the origin, then 
$$|I_{S,\phi}(\epsilon)| < C\epsilon^{1 \over d} \hskip 40 pt |J_{S,\phi}(\lambda)| < 
C|\lambda|^{-{1 \over d}} \eqno (4.12)$$

\noindent {\bf Proof.} Fix some $N$. By $(4.1a)-(4.1b)$, the upper bounds of $(4.12)$ hold for 
$I_{S,\phi}^{A_i}(\epsilon)$, $J_{S,\phi}^{A_i}(\lambda)$, $I_{S,\phi}^{B_i}(\epsilon)$, and  
$J_{S,\phi}^{B_i}(\lambda)$. Thus it suffices to prove these upper bounds for
$I_{S,\phi}^{D}(\epsilon)$ and $J_{S,\phi}^{D}(\epsilon)$, as well as $I_{S,\phi}^{C_i}(\lambda)$ and
$J_{S,\phi}^{C_i}(\lambda)$ if they derive from the vertical or horizontal ray. These are all done 
basically the same way, so we restrict our attention to $I_{S,\phi}^{D}(\epsilon)$ and 
$J_{S,\phi}^{D}(\epsilon)$. As in the proof of Lemma 4.1, we divide $D = D_1 \cup D_2$ along the curve
$y = x^m$. The two pieces are done similarly, so we will only consider $D_1$, the part where $y < x^m$.

Each vertical cross section $(D_1)_x$ of the set $D_1$ is the
a subset of the interval $[0,x^m]$. Hence there is some $\delta > 0$ such that on $D$, 
$\phi(x,y) = 0$ for $x < \delta$. Doing the $y$ integration first we have
$$|I_{S,\phi}^{D_1}(\epsilon)| < C\int_{\delta}^1 |\{y \in (D_1)_x: 0 < S(x,y) < \epsilon\}|\,dx \eqno (4.13)$$
By $(4.6)$, $|\partial_y^dS(x,y)|$ is bounded below on $x > \delta$. Hence by the Van der Corput
lemma in the $y$ direction, we have $|\{y \in (D_1)_x: 0 < S(x,y) < \epsilon\}| < C\epsilon^{1 \over d}$
uniformly in $x > \delta$. Inserting this back into $(4.13)$ gives the desired bounds. For the 
oscillatory integral, one similarly uses the Van der Corput lemma in the $y$ direction to get
$|\int_{(D_1)_x} e^{i\lambda S(x,y)}\phi(x,y)\,dy|  < C|\lambda|^{-{1 \over d}}$ uniformly in $x > \delta$.
Thus
$$|J_{S,\phi}^{D_1}(\lambda)| = |\int_\delta^1(\int_{(D_1)_x} e^{i\lambda S(x,y)}\phi(x,y)\,dy)\,dx| <
C|\lambda|^{-{1 \over d}} \eqno (4.14)$$
These are the sought-after bounds for $J_{S,\phi}^D(\lambda)$ and we are done.

We now proceed to the proof of the explicit formula $(1.7)$. The general 
methodology is similar to that of the case 1 arguments of section 3. If one writes $\phi =
\phi_1 + \phi_2$, where $\phi_1$ is supported in a smaller neighborhood of the origin and $\phi_2$ is 
zero on a neighborhood of the origin then by Lemma 4.3, $\limsup_{\epsilon \rightarrow 0} { I_{S,\phi}^+
(\epsilon) \over \ln(\epsilon)\epsilon^{1 \over d}} = \limsup_{\epsilon \rightarrow 0} { I_{S,\phi_1}^+
(\epsilon) \over \ln(\epsilon)\epsilon^{1 \over d}}$ and $\liminf_{\epsilon \rightarrow 0} { I_{S,\phi}^+
(\epsilon) \over \ln(\epsilon)\epsilon^{1 \over d}} = \liminf_{\epsilon \rightarrow 0} { I_{S,\phi_1}^+
(\epsilon) \over \ln(\epsilon)\epsilon^{1 \over d}}$. So when proving $(1.7)$ one can always replace 
$\phi$ by $\phi_1$ at will, regardless of how small the support of $\phi_1$ is. 

We will show if one first chooses the parameter $L$ of $(4.2b)$
sufficiently large, and then chooses $\phi_1$ to be supported in a sufficiently small neighborhood of the 
origin, then the above limsup and liminf, added to their analogues from the other three quadrants, can 
both be made as arbitrarily close to the limit given in 
$(1.7)$. To do this, in view of $(4.2b)$, it suffices to show that if the support of $\phi_1$ is 
sufficiently small, the quantities $\limsup_{\epsilon \rightarrow 0} { I_{S,\phi_1}^D \over \ln(\epsilon)
\epsilon^{1 \over d}}$ and $\liminf_{\epsilon \rightarrow 0} { I_{S,\phi_1}^D \over \ln(\epsilon)\epsilon^
{1 \over d}}$,  can be made arbitrarily close to the appropriate expression. For a fixed $L$ we will find lower 
bounds for the limsup and upper bounds for the liminf. In doing so, we will choose the parameter $M$
of the Taylor expansions in terms of $L$, and then the parameter $N$ of $(4.1a)-(4.1b)$ in terms of $L$ and
$M$. Analogous to in section 3 taking limits as $L$ goes to infinity, both expressions will converge to 
the same limit. Adding this limit to its analogues in the other 3 quadrants will give $(1.7)$.

We start with the following lemma. As before, $s_{dd} x^dy^d$ denotes the $(d,d)$ term of the 
Taylor expansion $S(x,y) = \sum_{a,b}s_{ab}x^ay^b$ of $S$ at the origin. 

\noindent {\bf Lemma 4.4.} There are constants $\beta$ and $C$
depending on $S(x,y)$, and a neighborhood $U$ of the origin depending on $S(x,y)$ and $L$,
such that $|S(x,y) - s_{dd}x^dy^d| < |C L^{-\beta}x^dy^d|$ on $D \cap U$. 

\noindent {\bf Proof.} Analogous to $(4.5)$, we may Taylor expand
$$S(x,y) -  s_{dd}x^dy^d = \sum_{M > a \geq d,\,\,\,M > b \geq d,\,\,\,(a,b) \neq (d,d)}s_{ab}'x^ay^b$$
$$ + \sum_{a < d,\,\,\,M > b > d,\,\,\,a + m_2b \geq \alpha_2}s_{ab}x^ay^b + \sum_{M > a > d,\,\,\,b < d,\,\,\,
a + m_1b \geq \alpha_1}s_{ab}x^ay^b + E_M(x,y) \eqno (4.15)$$
Here the second term is nonempty only if $(d,d)$ is the lower vertex of a compact edge, and the third 
term is nonempty only if $(d,d)$ is the upper vertex of a compact edge. 

The first sum can be made less than ${1 \over L}x^dy^d$ in absolute value by making the radius of $U$ sufficiently
small depending on $M$ and $S(x,y)$.
If the second sum is nonempty, then the domain $D$ is a subset of $\{(x,y) \in [0,1] \times [0,1]:0 < y <
N^{-m_2}x^{m_2}\}$. As a result if one changes coordinates from $(x,y)$ to $(x,y')$, where 
$y' = x^{m_2}y$, $D$ becomes a subset of $D' = \{(x,y) \in [0,1] \times [0,1]: 0 < y < N^{-m_2}\}$.
Observe that a given term $s_{ab}x^ay^b$ of the second sum becomes $s_{ab}x^{a + m_2b}(y')^b$. Since
$a + m_2b \geq \alpha_2$ and $b > d$ in each term in the second sum, the entire sum can be written as
$y'(x^{\alpha_2}
(y')^d)f(x,y')$ for some $f(x,y')$ which is a polynomial in $y$ and a fractional power of
$x$. Thus the sum is of absolute value at most $C_MN^{-m_2}x^{\alpha_2}(y')^d$ in a small enough 
neighborhood of the origin. Note that $s_{dd}x^dy^d = s_{dd}x^{d + dm_2}(y')^d$, and this is equal to 
$s_{dd}x^{\alpha_2}(y')^d$ since $(d,d)$ is on the edge with equation $t_1 + m_2t_2 = \alpha_2$. As a result,
in the original $(x,y)$ coordinates, the sum is of absolute value at most $C_MN^{-m_2}x^dy^d$. Thus 
if one chooses $N$ sufficiently large for fixed $L$ and $M$, one has the desired bounds.

The third sum is dealt with in exactly the same way, reversing the roles of the $x$ and $y$ axes. Since
$D$ necessarily lies in the range $x^{{1 \over L}} > y > x^L$, the error term $E_M(x,y)$ can be made
less than ${1 \over L} x^dy^d$ by making the radius of $U$ sufficiently small. This
completes the proof of Lemma 4.4.

\noindent {\bf Proof of $(1.7)$.} As before it suffices to assume $\phi(0,0) > 0$ as the general case 
can be obtained by writing $\phi = \phi' - \phi''$ where $\phi'(0,0)$ and  $\phi''(0,0) > 0$ are 
positive. Let $\delta_L > 0$ be such that on the ball $B(0,\delta_L)$ one has
$$(1 - L^{-\beta})\phi(0,0) < \phi(x,y) < (1 + L^{-\beta})\phi(0,0) \eqno (4.16)$$
Further assume that $\delta_L$ is small enough that $B(0,\delta_L) \subset U$, where $U$ is as in the 
previous lemma. Let $\psi(x,y)$ be a nonnegative cutoff function such that $0 \leq \psi \leq 1$, 
$\psi(x,y)$ is supported on $B(0,\delta_L)$, and $\psi(x,y) = 1$ on $[0,{\delta_L \over 2}] \times [0,{\delta_L \over 2}]$. Then by
the discussion following $(4.14)$ we may replace $\phi(x,y)$ by $\phi(x,y)\psi(x,y)$ without affecting 
the liminf or limsup. We have
$$I_{S,\phi}^D(\epsilon) < (1 + L^{-\beta})\phi(0,0) \big|\{(x,y) \in D: 0 < S(x,y) < \epsilon\}\big|$$
By Lemma 4.4 this is bounded by 
$$ (1 + L^{-\beta})\phi(0,0) \big|\{(x,y) \in D: 0 < (s_{dd} - CL^{-\beta})x^dy^d < \epsilon\}\big|
\eqno (4.17a)$$
In addition,
$$I_{S,\phi}^D(\epsilon) > (1 - L^{-\beta})\phi(0,0) \big|\{(x,y) \in D \cap ([0,{\delta_L \over 2}] 
\times [0,{\delta_L \over 2}]): 0 < S(x,y) < \epsilon\}\big|$$
$$ > (1 - L^{-\beta})\phi(0,0) \big|\{(x,y) \in D \cap [0,({\delta_L \over 2}] \times [0,{\delta_L 
\over 2}]):0 < (s_{dd} + CL^{-\beta})x^dy^d < \epsilon\}\big| \eqno (4.17b)$$
If $s_{dd} < 0$, then by $(4.17a)$ $I_{S,\phi}^D(\epsilon) = 0$ for large enough $L$. For $s_{dd} > 0$,
we need the following lemma, whose proof is routine.

\noindent {\bf Lemma 4.5}. Suppose $0 < m_2 < m_1$ and $\delta_0 > 0$. Then as $t \rightarrow 0$,
$$\big|\{(x,y) \in (0,\delta_0] \times (0,\delta_0]: x^{m_1} < y < x^{m_2},\,\,y < {t \over x}\}\big| = 
({1 \over m_1 + 1} - {1 \over m_2 + 1})t \ln(t) + O(t)$$
$$\big|\{(x,y) \in (0,\delta_0] \times (0,\delta_0]: 0 < y < x^{m_2},\,\,y < {t \over x}\}\big| = 
{-{1 \over m_2 + 1}}t \ln(t) + O(t)$$
We now apply Lemma 4.5 to $(4.17a)-(4.17b)$. We get that
$$I_{S,\phi}^D(\epsilon) < (1 + L^{-\beta})(s_{dd} - CL^{-\beta})^{-{1 \over d}}\phi(0,0)
\epsilon^{1 \over d}\ln(\epsilon)({1 \over m_1 + 1} - {1 \over m_2 + 1}) + C\epsilon^{1 \over d} \eqno (4.18a)$$
$$I_{S,\phi}^D(\epsilon) > (1 - L^{-\beta})(s_{dd} + CL^{-\beta})^{-{1 \over d}}\phi(0,0)
\epsilon^{1 \over d}\ln(\epsilon)({1 \over m_1 + 1} - {1 \over m_2 + 1}) - C\epsilon^{1 \over d}\eqno (4.18b)$$
(When $(d,d)$ is on the horizontal ray one substitutes $m_1 = L$, and when it is on the vertical ray one 
substitutes $m_2 = {1 \over L}$). Hence we have
$$\limsup_{\epsilon \rightarrow 0} {I_{S,\phi}^D(\epsilon) \over \epsilon^{1 \over d}\ln(\epsilon)}
\leq (1 + L^{-\beta})(s_{dd} - CL^{-\beta})^{-{1 \over d}}\phi(0,0)({1 \over m_1 + 1} - 
{1 \over m_2 + 1})\eqno (4.19a)$$
$$\liminf_{\epsilon \rightarrow 0} {I_{S,\phi}^D(\epsilon) \over \epsilon^{1 \over d}\ln(\epsilon)}
\geq (1 - L^{-\beta})(s_{dd} + CL^{-\beta})^{-{1 \over d}}\phi(0,0)({1 \over m_1 + 1} - 
{1 \over m_2 + 1}) \eqno (4.19b)$$
We take limits as $L \rightarrow \infty$. Both expressions converge to $\phi(0,0)({1 \over m_1 + 1} - 
{1 \over m_2 + 1})$, where now $m_1$ is taken as  $\infty$ when $(d,d)$ is on the horizontal ray and 
$m_2$ is taken as 0 when it is on the vertical ray. Hence by the discussion following $(4.14)$ we 
conclude that 
$$\lim_{\epsilon \rightarrow 0} {I_{S,\phi}^+(\epsilon) \over \epsilon^{1 \over d}\ln(\epsilon)}
= s_{dd}^{-{1 \over d}}\phi(0,0)({1 \over m_1 + 1} - {1 \over m_2 + 1})$$
Letting $s_i$ be the slope $-{1 \over m_i}$ of the edge $t_1 + m_it_2 = \alpha_i$, this becomes
$$\lim_{\epsilon \rightarrow 0} {I_{S,\phi}^+(\epsilon) \over \epsilon^{1 \over d}\ln(\epsilon)}
= s_{dd}^{-{1 \over d}}\phi(0,0)({1 \over s_1 - 1} - {1 \over s_2 - 1}) \eqno (4.20)$$
In summary, if $s_{dd} > 0$ then $(4.20)$ gives the contribution to 
$(1.7)$ from the upper right-hand quadrant. If $s_{dd} < 0$ then the contribution is zero as 
mentioned above Lemma 4.5. Adding
this to its analogues over the other three quadrants gives exactly the formula of $(1.7)$ and we are
done.

\noindent Our final task is to prove $(1.13b)$:

\noindent {\bf Proof of (1.13b).} Since the result is symmetric in $S$ and $-S$, we may replace
$S$ by $-S$ if necessary and assume that $s_{dd} > 0$. As in the proof of $(1.13a)$, we write $\phi = 
\phi_1 + \phi_2$, where $\phi_1$ is nonnegative and $\phi_2$ is zero on a neighborhood of the origin.
By Lemma 4.3, $\lim_{\epsilon \rightarrow 0} 
{I_{|S|,\phi_2}(\epsilon) \over \ln(\epsilon) \epsilon^{1 \over d}} = 0$. So to prove $(1.13b)$ it 
suffices to show $I_{|S|,\phi_1}(\epsilon) > A_{S,\phi_1}|\ln(\epsilon)|\epsilon^{1 \over d}$ for some 
$A_{S,\phi_1} > 0$. Since $I_{|S|,\phi_1}(\epsilon) \geq I_{S,\phi_1}(\epsilon)$, it further suffices 
to show that $I_{S,\phi_1}(\epsilon) > A_{S,\phi}|\ln(\epsilon)|\epsilon^{1 \over d}$ for some
$A_{S,\phi_1}$. For this we use Lemma 4.4, which implies that there is a 
$\delta > 0$ such that on $([0,\delta] \times [0,\delta]) \cap D$ we have
$$S(x,y) < {3 \over 2}s_{dd}x^dy^d \eqno (4.21)$$
As a result, shrinking $\delta$ to ensure that $\phi(x,y) > {1 \over 2} \phi(0,0)$ on $[0,\delta] 
\times [0,\delta]$ if necessary, we have 
$$I_{S,\phi_1}(\epsilon) > {1 \over 2}\phi(0,0) |\{(x,y) \in D \cap ([0,\delta] \times [0,\delta]):
{3 \over 2}s_{dd}x^dy^d < \epsilon\}|\eqno (4.22)$$
Using Lemma 4.5, we conclude that there is some $A_{S,\phi_1}$ with $|I_{S,\phi_1}(x,y)|> 
A_{S,\phi_1}|\ln(\epsilon)| \epsilon^{1 \over d}$ as needed. This gives $(1.13b)$ and we are done.

\noindent {\bf 5. Case 3 proofs.} 

In this section, $S(x,y)$ is a smooth phase function in case 3 of 
superadapted coordinates with $d > 1$. We restrict ourselves to the situation where the bisectrix 
intersects the horizontal ray in its interior, as the case of a vertical ray is entirely analogous.
Thus the lowest vertex of $N(S)$ is of the form $(c,d)$, where $c < d$ and $d$ is also the Newton distance
of $S$. As in sections 3 and 4 
we will focus our attention on the analysis of $I_{S,\phi}^+
(\epsilon)$ and $J_{S,\phi}^+(\lambda)$. We divide $[0,1] \times [0,1]$ into two parts. For a sufficiently
large positive integer $k$ 
(to be determined by our arguments), we let $D_1 = 
\{(x,y) \in [0,1] \times [0,1]: y < x^{k}\}$ and $D_2 = \{(x,y) \in [0,1] \times [0,1]: y > x^{k}\}$.
Our first lemma is the following.

\noindent {\bf Lemma 5.1.} $|I_{S,\phi}^{D_1}(\epsilon)| < C \epsilon^{1 \over d}$ and $|J_{S,\phi}^
{D_1}(\lambda)| < C|\lambda|^{-{1 \over d}}$.

\noindent {\bf Proof.} As in cases 1 and 2, we write the Taylor expansion of
$S$ at the origin as $S(x,y) = \sum_{a < M,\,\,\,b < M}s_{ab}x^ay^b + E_M(x,y)$, where for 
$0 \leq \alpha, \beta \leq M$ the function $E_M(x,y)$ satisfies the error estimates
$$|\partial_x^\alpha \partial_y^\beta E_M(x,y)| < C(|x|^{M - \alpha} + |y|^{M - \beta}) \eqno (5.1)$$
The $d$th $y$-derivative can be written as
$$\partial_y^d S(x,y) = \sum_{a < M,\,\,\,b < M - d}s_{ab}'x^ay^b + \partial_y^d E_M(x,y) \eqno (5.2)$$
Furthermore the Newton polygon of $\partial_y^d S(x,y)$ has a 
vertex at $(c,0)$, contained either in the vertical ray of $N(\partial_y^d S)$ or an edge of 
$N(\partial_y^d S)$ with equation $t_1 + mt_2 = c$ with $m > 0$. Hence each $(a,b)$ in the sum of $(5.2)$
satisfies $a + mb \geq c$. Analogous to $(4.5)$ we rewrite $(5.2)$ as
$$\partial_y^d S(x,y) = s_{c0}'x^c + \sum_{a < c,\,\,\,0 < b < M - d,\,\,\,a + mb \geq c}s_{ab}'x^ay^b $$
$$+ \sum_{c \leq a < M,\,\,\,0 \leq b < M - d,\,\,\,(a,b) \neq (c,0)}s_{ab}'x^ay^b + 
{\partial_y^d}E_M(x,y) \eqno (5.3)$$
In the case where $(c,0)$ is on a vertical ray of $N(\partial_y^d S)$, the first sum of $(5.3)$ is empty. We now argue 
like after $(4.5)$. Since $y < x^{k}$ for all $(x,y) \in D_1$, if $(x,y)$ is in $D_1$ each term
$s_{ab}'x^ay^b$ in the first sum of $(4.5)$ is bounded in absolute value by $|s_{ab}'|x^{a + kb} 
\leq |s_{ab}'|x^{k-m}(x^{a + mb}) \leq |s_{ab}'|x^{k - m + c}$.  Thus as long as $k$ were chosen 
greater than $m$, which we may assume, then if the support of
$\phi(x,y)$ is sufficiently small, then for $(x,y)$ in this support, the absolute value of the whole
first sum is less than ${1 \over 4}|s_{c0}'|x^c$. Also, the absolute value of a given term 
$s_{ab}'x^ay^b$ in the second sum is at most $|s_{ab}'|x^c(x + y)$. As a result, for such $(x,y)$ the
absolute value of the second sum is also at most ${1 \over 4}|s_{c0}'|x^c$. Similarly, using $(5.1)$
and the fact that $0 < y < x^{k}$, for such $(x,y)$ the quantity $|{\partial_y^d}E_M(x,y)|$ can also be 
assumed to be at most ${1 \over 4}|s_{c0}'|x^c$. Consequently, for these $(x,y)$ we have
$$|\partial_y^d S(x,y)| > {1 \over 4}|s_{c0}'|x^c \eqno (5.4)$$
As a result, by the measure version of Van der Corput's lemma of [C], for each $x$ in this range we have
$$|\{y : 0 < S(x,y) < \epsilon\}| < C\epsilon^{1 \over d}x^{-{c \over d}} \eqno (5.5)$$
Consequently, integrating with respect to $y$ first one has
$$|I_{S,\phi}^{D_1}(\epsilon)| < C\int_0^1 \epsilon^{1 \over d}x^{-{c \over d}}\,dx = 
C' \epsilon^{1 \over d} \eqno (5.6)$$
This is the desired upper bound for $I_{S,\phi}^{D_1}(\epsilon)$. As for the oscillatory integral 
analogue, by $(5.4)$ the normal Van der Corput lemma in the $y$ direction gives
$$|\int e^{i \lambda S(x,y)} \phi(x,y) \,dy| \leq C |\lambda|^{-{1 \over d}}x^{-{c \over d}}\eqno (5.7)$$
Hence by integrating first with respect to $y$ one has 
$$|J_{S,\phi}^{D_1}(\lambda)| < C\int_0^1 |\lambda|^{-{1 \over d}}x^{-{c \over d}}\, dx = 
C' |\lambda|^{-{1 \over d}} \eqno (5.8)$$
This completes the proof of Lemma 5.1.

\noindent It turns out that one gets stronger estimates for the $I_{S,\phi}^{D_2}(\epsilon)$ and 
$J_{S,\phi}^{D_2}(\lambda)$. Observe that since $c < d$, the quantity ${1 + k \over c + kd}$ is greater
than ${1 \over d}$. We have the following.

\noindent {\bf Lemma 5.2.} $|I_{S,\phi}^{D_2}(\epsilon)| < C \epsilon^{1 + k \over c + kd}$ and 
$|J_{S,\phi}^{D_2}(\lambda)| < C|\lambda|^{-{1 + k \over c + kd}}$.

\noindent {\bf Proof.} We will verify the hypotheses of Lemma 2.0, with the roles of the $x$ and the
$y$ variables reversed. Because $(c,d)$ is the rightmost vertex of $N(S)$, the Newton polygon of 
$\partial_x^cS(x,y)$ has a single vertex at $(0,d)$. The Taylor expansion of 
$\partial_x^cS(x,y)$ can be written in the form
$$\partial_x^cS(x,y) = r_{0d}\,y^d + \sum_{0 \leq a < M,\,\,\,d \leq b < M\,\,\, (a,b) \neq (0,d)}
r_{ab}x^ay^b + E_M(x,y) \eqno (5.9)$$
Similar to elsewhere in this paper, bounding the error term using the fact that $x < y^{{1 \over k}}$
on $D_2$, in a small enough neighborhood of the origin on $D_2$ one has
$$|\partial_x^cS(x,y)| > {1 \over 2}|r_{0d}|y^d $$
Thus if $c \geq 2$, one can apply Lemma 2.0 and immediately get this lemma. If $c = 1$, to apply 
Lemma 2.0 one also needs
that $|\partial_x^2S(x,y)| <  C y^{d - {1 \over k}}$. But in fact since the Newton polygon of 
$\partial_x^2S(x,y)$ is a subset of $\{(x,y): y \geq d\}$, by expanding as in $(5.9)$ one even has 
the stronger estimate 
$|\partial_x^2S(x,y)| <  C y^d$. Thus Lemma 2.0 applies here. If $c = 0$, to apply Lemma 2.0 one needs
$(2.2c)$ to hold (with the $x$ and $y$ variables reversed) which here means one needs
$|\partial_y S(x,y)| >  C y^{d - 1}$ and  $|\partial_y^2 S(x,y)| <  C' y^{d - 2}$. Since
the Newton polygon of $\partial_y^i S(x,y)$ has a single vertex at $(0,d-i)$ this holds as in the 
$c = 1$ case. Lastly, to apply Lemma 2.0  for $c = 0$ one also needs that $d > {1 \over k} + 1$. We 
can make this true simply choosing $k$ large enough since $d$ is at least 2. This completes the proof
of Lemma 5.2.

\noindent {\bf Corollary 5.3.} $(1.11)$ and $(1.15a)$ hold in Case 3.

\noindent {\bf Proof.} Add the estimates from Lemmas 5.1 and 5.2 and their analogues from the other 
three quadrants.

Next, we prove the lower bounds of $(1.13c)$. Since we are not trying to prove sharp estimates, the 
arguments are not that intricate. Assume $\phi(0,0) \neq 0$, and let $M$ be some large positive integer.
We examine the behavior of $S(x,y)$ on the set $Z_N = \{(x,y): x > 0,\,\,x^N < y < 2x^N \}$ for $N$ 
sufficiently large. We Taylor expand $S(x,y) = \sum_{a < M,\,\,\,b < M}s_{ab}x^ay^b + E_M(x,y)$ as above.
If $N$ is large enough, the term $s_{cd}x^cy^d$ dominates this Taylor expansion much the way $r_{0d}y^d$
dominates $(5.9)$ or $s_{c0}'x^c$ dominates $(5.3)$. Hence in a small enough neighborhood $U$ of the 
origin, on $Z_N$ we have
$$|S(x,y)| < 2|s_{cd}|x^cy^d \eqno (5.10)$$
Shrinking $U$, we may assume that $|\phi(x,y)| > {1 \over 2}|\phi(0,0)|$. Hence for
$c_{\phi} = {1 \over 2}|\phi(0,0)|$ we have
$$I_{|S|,\phi}(\epsilon) > c_{\phi} |\{(x,y) \in U \cap Z_N: 2|s_{cd}|x^cy^d < \epsilon \}| \eqno (5.11)$$
It is easy to compute that the curve $y = 2x^N$ intersects the curve $2|s_{cd}|x^cy^d = \epsilon$ at 
$x = a\epsilon^{1 \over Nd + c}$ for some $a$ depending on $S(x,y)$. Hence for $\epsilon$ sufficiently 
small the measure of the set in
the right-hand side of $(5.11)$ is at least the measure of the portion of $Z_N$ between $x = {a \over 2}
\epsilon^{1 \over Nd + c}$ and $x = a\epsilon^{1 \over Nd + c}$, given by $a'\epsilon^{N + 1 \over
Nd + c}$ where now $a'$ also depends on $N$. Thus we can write
$$I_{|S|,\phi}(\epsilon) > c'_{S,\phi, N} \epsilon^{N + 1 \over Nd + c} \eqno (5.12)$$
Since $d > c$, the exponent in $(5.12)$ is larger than ${1 \over d}$ but as $N \rightarrow \infty$ it
tends to ${1 \over d}$. This gives us $(1.13c)$.

We now move to the case of real-analytic phase. Our goal here is to prove Theorem 1.1 c).
So assume $S(x,y)$ is real-analytic. It suffices to show that $\lim_{\epsilon 
\rightarrow 0} {I_{S,\phi}(\epsilon) \over \epsilon^{1 \over d}}$ exists and is given by $(1.8)$.
As in cases 1 and 2, we will give an expression for $\lim_{\epsilon 
\rightarrow 0} {I_{S,\phi}^+(\epsilon) \over \epsilon^{1 \over d}}$ and the full limit will follow by 
adding this and the analogues from the other quadrants. Also, by Lemma 5.2, $\lim_{\epsilon 
\rightarrow 0} {I_{S,\phi}^{D_2}(\epsilon) \over \epsilon^{1 \over d}} = 0$ so it suffices to show
$\lim_{\epsilon \rightarrow 0} {I_{S,\phi}^{D_1}(\epsilon) \over \epsilon^{1 \over d}}$ exists and has
the desired value.
Because the bisectrix intersects $N(S)$ in the interior of its horizontal ray and $(c,d)$ is the 
lowest vertex of $N(S)$, the real-analytic $S(x,y)$ can be written as 
$$S(x,y) = s_{cd}x^cy^d + x^{c+1}y^dg(x) + \sum_{b \geq d+1} s_{ab}x^ay^b \eqno (5.13a)$$
Here $g(x)$ is real-analytic. Changing coordinates from $(x,y)$ to $(x,y')$ where $y = x^ky'$, we 
have 
$$S(x,x^{k}y') = s_{cd}x^{c + kd}(y')^d + x^{c +  kd + 1}(y')^dg(x)  +  \sum_{b \geq d+1} 
s_{ab} x^{a + kb}(y')^b \eqno (5.13b)$$
Since the line $ t_1 + mt_2 = \alpha$ is an edge of $N(S)$ containing $(c,d)$, each $(a,b)$ in the sum
$(5.13a)$ satisfies $a + mb \geq c + md$. Furthermore, $b > d$, and therefore
$$a + kb = a + mb + (k-m)b \geq c + md + (k-m)b > c + md + (k-m)d = c + kd$$
Hence we can rewrite $\sum_{b \geq d+1} s_{ab} x^{a + kb}(y')^b$ as $x^{c + kd + 1}
(y')^{d + 1}f(x,y)$ where $f(x,y)$ is real-analytic. Thus we have
$$S(x,x^{k}y') = s_{cd}x^{c + kd}(y')^d + x^{c +  kd + 1}(y')^dg(x)  +  x^{c + kd + 1}(y')^{d + 1}
f(x,y') \eqno (5.14)$$
In the $(x,y')$ coordinates $I_{S,\phi}^{D_1}(\epsilon)$ becomes
$$I_{S,\phi}^{D_1}(\epsilon) = \int_{\{(x,y') \in [0,1] \times [0,1]: 0 < S(x,x^{k}y') < \epsilon\}} 
x^{k} \phi(x,x^{k}y')\,dx\,dy' \eqno (5.15)$$
Let $d'$ be between $c$ and $d$. The exact value of $d'$ will be dictated by our arguments.
Then the portion of $(5.15)$ over $x < \epsilon^{1 \over d'(k + 1)}$
has absolute value at most $C\int_0^{\epsilon^{1 \over d'(k + 1)}} x^{k}\,dx = 
C\epsilon^{1 \over d'}$.
Since this is $o(\epsilon^{1 \over d})$, this portion of the integral will can be removed without
affecting $\lim_{\epsilon \rightarrow 0} {I_{S,\phi}^{D_1}(\epsilon) \over \epsilon^{1 \over d}}$. 
In other words, we may replace $I_{S,\phi}^{D_1}(\epsilon)$ by $I_{S,\phi}'(\epsilon)$ where
$$I_{S,\phi}'(\epsilon) = \int_{\{(x,y') \in [\epsilon^{1 \over d'(k + 1)},1] \times [0,1] : 0 < 
S(x,x^{k}y') < \epsilon\}} x^{k} \phi(x,x^{k}y')\,dx\,dy' \eqno (5.16)$$
We now fix $x > \epsilon^{1 \over d'(k + 1)}$ and look at the set $E_x = \{y \in [0,1]: 0
< S(x,x^{k}y') < \epsilon\}$. We may assume the support of $\phi(x,y)$ is small enough so that if
$\phi(x,y) \neq 0$ then $x$ is small enough so that 
$$|x^{c +  kd + 1}(y')^dg(x)| +  |x^{c + kd + 1}(y')^{d + 1}f(x,y)| < {1 \over 2}|s_{cd}|x^{c + kd}y^d
\eqno (5.17)$$
Note that if $s_{cd}$ is negative, by $(5.14)$ and $(5.17)$  $S(x,x^{k}y')$ is always
negative and thus $I_{S,\phi}'(\epsilon) = 0$. So assume that $s_{cd} > 0$.
Then by $(5.14)$ and $(5.17)$ we have
$$E_x \subset \{y \in [0,1]: 0 <  s_{cd}x^{c + kd}y^d < 2\epsilon\} \subset [0, C\epsilon^
{1 \over d}x^{-{c + kd \over d}}]\eqno (5.18)$$
Since $x \geq \epsilon^{1 \over d'(k + 1)}$, we have
$$\epsilon^{1 \over d}x^{-{c + kd \over d}} \leq \epsilon^{1 \over d}\epsilon^{-{c + kd \over d(d')
(k + 1)}} = \epsilon^{kd' + d' \over d(d')(k + 1)}\epsilon^{-{c + kd \over d(d')(k + 1)}} = 
\epsilon^{k(d' - d) + (d' - c) \over dd'(k + 1)}$$
Hence if $d'$ were chosen close enough to $d$, there is some $\eta' > 0$ such that for $x \geq 
\epsilon^{1 \over d'(k + 1)}$ one has $\epsilon^{1 \over d}x^{-{c + kd \over d}} < \epsilon^
{\eta'}$ and thus
$$E_x \subset [0,C\epsilon^{\eta'}] \eqno (5.19)$$
Next, we write $I_{S,\phi}'(\epsilon) = I_{S,\phi}''(\epsilon) + I_{S,\phi}'''(\epsilon)$, where
$$I_{S,\phi}''(\epsilon) = \int_{\{(x,y') \in [\epsilon^{1 \over d'(k + 1)},1] \times [0,1] : 0 < 
S(x,x^{k}y') < \epsilon\}} x^{k} \phi(x,0)\,dx\,dy' \eqno (5.20a)$$
$$I_{S,\phi}'''(\epsilon) = \int_{\{(x,y') \in [\epsilon^{1 \over d'(k + 1)},1] \times [0,1] : 0 < 
S(x,x^{k}y') < \epsilon\}} x^{k} (\phi(x,x^{k}y') - \phi(x,0))\,dx\,dy' \eqno (5.20b)$$
Note that due to $(5.19)$, the factor $(\phi(x,x^{k}y') - \phi(x,0))$ in $(5.20b)$ is bounded in
absolute value by
$C'\epsilon^{\eta'}$, so we have
$$I_{S,\phi}'''(\epsilon) \leq C'\epsilon^{\eta'}\int_{\{(x,y') \in [\epsilon^{1 \over d'(k + 1)},1]
\times [0,1] : 0 < S(x,x^{k}y') < \epsilon\}} x^{k} \,dx\,dy' \eqno (5.21)$$
In $(5.21)$, we perform the $y$ integration by inserting the second inclusion of $(5.18)$. We then have
$$I_{S,\phi}'''(\epsilon) \leq C''\epsilon^{\eta'}\int_{\epsilon^{1 \over d'(k + 1)}}^1
x^{k} \epsilon^{1 \over d}x^{-{c + kd \over d}}\,dx\,dy' \leq 
C'''\epsilon^{{1 \over d} + \eta'} \int_0^1x^{-{c \over d}}\,dx = C''''\epsilon^{{1 \over d} + \eta'}
\eqno (5.22)$$
Thus $\lim_{\epsilon \rightarrow 0} {I_{S,\phi}'''(\epsilon) \over \epsilon^{1 \over d}} = 0$. 
Hence $\lim_{\epsilon \rightarrow 0} {I_{S,\phi}'(\epsilon) \over \epsilon^{1 \over d}}
= \lim_{\epsilon \rightarrow 0} {I_{S,\phi}''(\epsilon) \over \epsilon^{1 \over d}}$, and our goal
now becomes to prove the latter limit gives the portion of $(1.8)$ coming from
the upper right quadrant. Next, we rewrite $I_{S,\phi}''(\epsilon)$ as
$$I_{S,\phi}''(\epsilon) = \int_{\epsilon^{1 \over d'(k + 1)}}^1 x^{k} \phi(x,0)|\{y': 0 < 
S(x,x^{k}y') < \epsilon\}| \,dx\,dy' \eqno (5.23)$$
To analyze $(5.23)$, first note that ${x^{c + kd + 1} \over s_{cd}x^{c + kd} +
x^{c +  kd + 1}g(x)}$ is a real-analytic function in a neighborhood of the origin, which we 
denote by $h(x)$. Then by $(5.14)$ we have
$$S(x,x^{k}y') = (s_{cd}x^{c + kd} +  x^{c +  kd + 1}g(x))((y')^d   +  xh(x)(y')^{d+1}f(x,y'))
\eqno (5.24)$$
Since $s_{cd}$ is being assumed to be positive, in $(5.24)$ $(s_{cd}x^{c + kd} +  x^{c +  kd + 1}g(x))^
{1 \over d}$ is positive and there is $j(x,y')$ with $\partial_yj(x,0) = 1$ and $\partial_xj(x,0) = 0$
such that $(5.24)$ can be rewritten as
$$S(x,x^{k}y')^{1 \over d} = (s_{cd}x^{c + kd}+  x^{c +  kd + 1}g(x))^{1 \over d}j(x,y') \eqno (5.25)$$
Consequently, in $(5.23)$, one has
$$|\{y': 0 < S(x,x^{k}y') < \epsilon\}| = \big|\{y': 0 < j(x,y') < 
\epsilon^{1 \over d}(s_{cd}x^{c + kd}+  x^{c +  kd + 1}g(x))^{-{1 \over d}}\}\big| \eqno (5.26)$$
By the inverse function theorem, $(x,j(x,y'))$ has an inverse function which
can be written as $(x, k(x,y'))$ for some $k(x,y')$ which satisfies $\partial_yk(x,0) = 1$ and 
$\partial_xk(x,0) = 0$ By $(5.19)$ the interval of $(5.26)$ has length at most $C\epsilon^{\eta'}$.
As a result, as long as $\epsilon$ is small enough we can use a linear approximation to $k(x,y')$ and
get that 
$$\big|\{y': 0 < j(x,y') < \epsilon^{1 \over d}(s_{cd}x^{c + kd}+  x^{c +  kd + 1}
g(x))^{-{1 \over d}}\}\big|$$
$$= \epsilon^{1 \over d}(s_{cd}x^{c + kd}+  x^{c +  kd + 1}g(x))^{-{1 \over d}} + \epsilon^{\eta'}
O(\epsilon^{1 \over d}(s_{cd}x^{c + kd}+  x^{c +  kd + 1}g(x))^{-{1 \over d}}) \eqno (5.27)$$
Thus we have
$$I_{S,\phi}''(\epsilon) = \int_{\epsilon^{1 \over d'(k + 1)}}^1 x^{k} \phi(x,0)\epsilon^
{{1 \over d}}(s_{cd}x^{c + kd} + x^{c +  kd + 1}g(x))^{-{1 \over d}} \,dx$$
$$+ O(\epsilon^{\eta'}\int_{\epsilon^{1 \over d'(k + 1)}}^1 x^{k}|\phi(x,0)|\epsilon^
{{1 \over d}}(s_{cd}x^{c + kd} + x^{c +  kd + 1}g(x))^{-{1 \over d}}\,dx) \eqno (5.28)$$
Because $\eta' > 0$, the ratio of the second term to $\epsilon^{{1 \over d}}$ goes to zero as 
$\epsilon > 0$ (assuming the integral is finite, which we will see shortly). Thus the second term
does not contribute to $\lim_{\epsilon \rightarrow 0} {I_{S,\phi}''(\epsilon) \over \epsilon^{1
\over d}}$ and we have
$$\lim_{\epsilon \rightarrow 0} {I_{S,\phi}''(\epsilon) \over \epsilon^{1 \over d}} = \int_0^1
x^{k}(s_{cd}x^{c + kd} + x^{c +  kd + 1}g(x))^{-{1 \over d}} \phi(x,0)\,dx $$
$$= \int_0^1 (s_{cd}x^c + x^{c + 1}g(x))^{-{1 \over d}}\phi(x,0) \,dx \eqno (5.29)$$
Since $|s_{cd}x^{c + kd} + x^{c +  kd + 1}g(x)| > {1 \over 2}x^{c + kd}$ and $c < d$ the integral
$(5.29)$ is finite as needed.
Going back to the definition of $g(x)$, $S(x,y) = (s_{cd}x^c + x^{c+1}g(x))
y^d +O(y^{d+1})$. As a result, $(5.29)$ translates into the part of equation $(1.8)$ coming from
the upper right-hand quadrant. (Recall that $(5.29)$ is for $s_{cd} > 0$ and that the limit
is zero when $s_{cd} < 0$). Adding the analogous expressions from the remaining three quadrants gives
$(1.8)$ and we are done.

\noindent {\bf 6. Proof of Theorem 1.6b.} 

Suppose $S(x,y)$ is a smooth phase function in case 1 superadapted coordinates, and $\phi(x,y)$ is 
nonnegative with $\phi(0,0) > 0$. Let $\psi(t)$ be a nonnegative function in $C_c(\R)$ 
such that $\psi(t) > 1$ on $[-1,1]$. Then by $(1.13a)$, if $j$ is sufficiently large we have
$$|\int \psi(2^jS(x,y))\phi(x,y)\,dx\,dy| > |\int_{\{(x,y): |S(x,y)| < 2^{-j}\}} \phi(x,y) \,dx\,dy|
> A_{S,\phi}2^{-{j \over d}} \eqno (6.1)$$ 
We also have
$$\int \psi(2^jS(x,y))\phi(x,y)\,dx\,dy = 2^{-j}\int(\int\hat{\psi}(2^{-j}\lambda)e^{i\lambda S(x,y)} 
d\lambda)dx\,dy$$
$$= 2^{-j}\int\hat{\psi}(2^{-j}\lambda)(\int e^{i\lambda S(x,y)} \phi(x,y)\,dx\,dy)d\lambda
= 2^{-j}\int\hat{\psi}(2^{-j}\lambda)J_{S,\phi}(\lambda)\,d\lambda \eqno (6.2)$$
In order to prove $(1.16a)$, we argue by contradiction. Suppose that we were in the setup of $(1.16a)$
but
$\limsup_{\lambda \rightarrow \infty}\big|{J_{S,\phi}(\lambda) \over \lambda^{-{1 \over d}}}\big| = 
0$. Then for any $\delta > 0$, we may let $M_{\delta}$ be such that for $|\lambda| > M_{\delta}$ we have
$|J_{S,\phi}(\lambda)| < \delta |\lambda|^{-{1 \over d}}$. We then have
$$|2^{-j}\int\hat{\psi}(2^{-j}\lambda)J_{S,\phi}(\lambda)\,d\lambda| = $$
$$|2^{-j}\int_{|\lambda| < M_{\delta}}\hat{\psi}(2^{-j}\lambda)J_{S,\phi}(\lambda)\,d\lambda + 
2^{-j}\int_{|\lambda| > M_{\delta}}\hat{\psi}(2^{-j}\lambda)J_{S,\phi}(\lambda)\,d\lambda|$$
$$< C 2^{-j}M_{\delta} + 2^{-j}\delta \int_{|\lambda| > M_{\delta}}|\hat{\psi}(2^{-j}\lambda)|
|\lambda|^{-{1 \over d}}\,d\lambda\eqno (6.3)$$
In turn, equation $(6.3)$ is bounded by 
$$C 2^{-j}M_{\delta} + 2^{-j}\delta \int_{\R}|\hat{\psi}(2^{-j}\lambda)||\lambda|^{-{1 \over d}}
\,d\lambda$$
$$= C 2^{-j}M_{\delta} + 2^{-{j \over d}}\delta \int_{\R}|\hat{\psi}(\lambda)||\lambda|^{-{1 \over d}}
d\lambda$$
$$ < C 2^{-j}M_{\delta} + C'\delta 2^{-{j \over d}} \eqno (6.4)$$
When $j$ is sufficiently large, $(6.4)$ is at most $2C'\delta 2^{-{j \over d}}$. On the other hand,
by $(6.1)$, it must also be at least $A_{S,\phi}2^{-{j \over d}}$. This gives a contradiction if 
$\delta$ were chosen less than ${A_{S,\phi} \over 2C'}$. This contradiction implies
that $\limsup_{\lambda \rightarrow \infty}\big|{J_{S,\phi}(\lambda) \over \lambda^{-{1 \over d}}}\big|$
is in fact positive, giving $(1.16a)$.

That the $\limsup$ of $(1.16c)$ is positive is proven from $(1.13c)$ exactly as $(1.16a)$ is proven from 
$(1.13a)$, so we do not include a proof here. Since it holds for all $\delta > 0$ the $\limsup$ is
automatically infinite. Equation $(1.16b)$ is proved similarly to $(1.16a)$, using $(1.13b)$ in place
of $(1.13a)$. Namely, suppose $S(x,y)$ is a smooth phase function in case 2
superadapted coordinates, and $\phi(x,y)$ is a nonnegative function with $\phi(0,0) > 0$.
Using $(1.13b)$ we have 
$$|\int \psi(2^jS(x,y))\phi(x,y)\,dx\,dy| > |\int_{\{(x,y): |S(x,y)| < 2^{-j}\}} \phi(x,y) \,dx\,dy|
> A_{S,\phi}{j \over d}2^{-{j \over d}} \eqno (6.5)$$
Exactly as above we also have 
$$\int \psi(2^jS(x,y))\phi(x,y)\,dx\,dy = 2^{-j}\int\hat{\psi}(2^{-j}\lambda)J_{S,\phi}(\lambda)\,
d\lambda \eqno (6.6)$$
Proceeding by contradiction again, suppose $(1.16b)$ does not hold. 
Therefore for every $\delta > 0$ there is some $L_{\delta}$ such that for $|\lambda| > L_{\delta}$ we have
$|J_{S,\phi}(\lambda)|<\delta |\lambda|^{-{1 \over d}}\ln|\lambda| $. Analogous to $(6.3)$ we have
$$2^{-j}|\int\hat{\psi}(2^{-j}\lambda)J_{S,\phi}(\lambda)\,d\lambda|\,\,\,<\,\,\,C 2^{-j}L_{\delta} + 
2^{-j}\delta \int_{|\lambda| > L_{\delta}}|\hat{\psi}(2^{-j}\lambda)||\lambda|^{-{1 \over d}}\ln|\lambda|
\,d\lambda$$
$$\leq C 2^{-j}L_{\delta} + 2^{-j}\delta \int_{\R}|\hat{\psi}(2^{-j}\lambda)| |\lambda|^{-{1 \over d}}\ln|\lambda|
\,d\lambda \eqno (6.7)$$
Changing variables, this in turn is equal to 
$$ C 2^{-j}L_{\delta} + 2^{-{j \over d}}\delta \int_{\R}|\hat{\psi}(\lambda)||\lambda|^{-{1 \over d}}
\ln(2^j|\lambda|)d\lambda$$
$$= C 2^{-j}L_{\delta} +  2^{-{j \over d}}\delta \int_{\R}|\hat{\psi}(\lambda)||\lambda|^{-{1 \over d}}
\ln|\lambda|\,d\lambda + j2^{-{j \over d}}\ln(2)\delta \int_{\R}|\hat{\psi}(\lambda)|
|\lambda|^{-{1 \over d}} \,d\lambda$$
$$< C 2^{-j}L_{\delta} + C'' \delta j 2^{-{j \over d}} \eqno (6.8)$$
If $j$ is sufficiently large, $(6.4)$ is at most $2C''\delta j 2^{-{j \over d}}$, while by 
$(6.5)$ it is at least $A_{S,\phi}{j \over d}2^{-{j \over d}}$. This is a contradiction if 
$\delta < {A_{S,\phi} \over 2 d C''}$. Therefore $\limsup_{\lambda \rightarrow 
\infty} \big|{J_{S,\phi}(\lambda) \over \ln(\lambda)\lambda^{-{1 \over d}}}\big|$ must in fact be 
positive. Hence we have $(1.16b)$ and we are done.

\noindent {\bf 7. Superadapted coordinates.} 

In this section we prove the existence of superadapted
coordinates for smooth phase functions. Here we always assume $S(x,y)$ is a smooth phase function defined on a
neighborhood of the origin such that $S(0,0) = 0$ and $S(x,y)$ has nonvanishing Taylor expansion at
the origin. The three cases of superadapted coordinates can be written as follows:

\noindent {\bf Case 1.} The bisectrix intersects $N(S)$ in the interior of a bounded edge 
$e$ and any real zero $r \neq 0$ of $S_e(1,y)$ or $S_e(-1,y)$ has order less than $d(S)$.

\noindent {\bf Case 2.} The bisectrix intersects $N(S)$ at a vertex $(d,d)$ and if $e$ is a compact edge
of $N(S)$ containing $(d,d)$ then any real zero $r \neq 0$ of $S_e(1,y)$ or $S_e(-1,y)$ has 
order less than $d(S)$. 

\noindent {\bf Case 3.} The bisectrix intersects $N(S)$ in the interior of one of the unbounded
edges.

\noindent {\bf Lemma 7.0.} Any superadapted coordinate system is adapted.

\noindent {\bf Proof.} By the main theorem of [G1], if $U$ is a small enough neighborhood of the origin
and $\epsilon_0$ denotes the supremum
of the numbers $\epsilon$ for which $\int_U |S|^{-\epsilon}$ is finite, then $d(S) \leq {1 \over 
\epsilon_0}$, with $d(S) = {1 \over \epsilon_0}$ in cases 1, 2, and 3. Hence if one is in cases 1, 2, 
or 3, one is in adapted coordinates.

\noindent Case 2 has some special features for which the following preliminary lemma will be useful. Related
lemmas occur in [PSSt] and [V]. 

\noindent {\bf Lemma 7.1.} Suppose the bisectrix intersects $N(S)$ at a vertex $(d,d)$ but is not
in superadapted coordinates. Correspondingly, let $e$ be a compact edge of $N(S)$ containing $(d,d)$ such
that $S_e(1,y)$ or $S_e(-1,y)$ has a zero of order $d$ or greater. If $(d,d)$ is the upper
vertex of $e$, then $S_e(x,y)$ is of the form 
$cx^{\alpha}({y \over x^m} - r)^d$ for positive integers $\alpha, m$ and some nonzero $c, r$. If $(d,d)$ 
is the lower vertex of $e$, then $S_e(x,y)$ has the analogous form $c'y^{\alpha'}({x \over y^{m'}} - r')^d$.
with $c',r' \neq 0$ and $\alpha', m'$ positive integers.

\noindent {\bf Proof.} We first consider the case where $(d,d)$ is the upper vertex of $e$. 
Write the equation of $e$ as $t_1 + mt_2 = \alpha$. We will show that these values
of $m$ and $\alpha$ work. Note that if $s_{ab}x^ay^b$ appears in $S_e(x,y)$ then $a + mb = \alpha$.
We factor out $x^{\alpha}$, writing $S_e(x,y) = x^{\alpha}T_e(x,y)$. Each term of $T_e(x,y)$ is
now of the form $t_{ab}x^{a - \alpha}y^b$ with  $(a - \alpha) + mb = 0$ or $(a - \alpha) = -mb$.
Thus we have
$$t_{ab}x^{a - \alpha}y^b = t_{ab}({y \over x^m})^b \eqno (7.1)$$
Conequently for a polynomial $P(z)$, we can write
$$S_e(x,y)= x^{\alpha}P({y \over x^m}) \eqno (7.2)$$
Plugging in $x = 1$ or $-1$, we see that $P(y)$ is a polynomial of degree $d$ with a real zero $r$ of order 
$d$ or 
greater. Therefore we must have $P(y) = c(y - r)^d$ for some $c \neq 0$. Hence  $S_e(x,y) = 
cx^{\alpha}({y \over x^m} - r)^d$. Since $x$ can only appear to integer powers, $m$ must be an integer
and therefore $\alpha$ is as well.
Also, since it comes from an edge $S_e(x,y)$ contains multiple terms. Hence $r \neq 0$ and we are
done with the case where $(d,d)$ is the upper vertex of $e$.

The case where $(d,d)$ is the upper vertex of $e$ is done similarly. Since $S_e(1,y)$ or $S_e(-1,y)$ has 
a zero $r \neq 0$
of order $d$ or more, $S_e(x,y)$ has zeroes of order $d$ along some curve $y = rx^m$. Hence $S_e(x,1)$ or
$S_e(-x,1)$ has a zero not at the origin of order $d$ or more and the above argument applies, reversing
the roles of the $x$ and $y$ variables.

\noindent {\bf Lemma 7.2.} Suppose one is not in superadapted coordinates. Suppose $e$ is an edge of
$N(S)$ intersecting the bisectrix in its interior with equation $t_1 + mt_2 = \alpha$ for $m \geq 1$
such that $S_e(1,y)$ has a zero $r \neq 0$ of order $k \geq d(S)$. Then $m$ is an integer and both 
$S_e(1,y)$ and $S_e(-1,y)$ have a zero of order $k$ not at the origin. 

\noindent {\bf Proof.}   
Exactly as $(7.2)$, there is some polynomial $Q(y)$ such that for $x > 0$ we have
$$S_e(x,y) = x^{\alpha}Q({y \over x^m}) \eqno (7.2')$$
Plugging in $x = 1$, we see that $Q(y) = S_e(1,y)$. 

We now show that $m$ must in fact be an integer. To see this, note that if $m$ were not an integer, 
then the degrees of the powers of $y$
appearing in $S_e(1,y)$ would have to be separated by at least 2. Hence $S_e(1,y)$ would have to be 
of the form $y^{\beta}R(y^c)$ for some $\beta \geq 0,$ $c \geq 2$, where $R$ is a polynomial.  Next, since 
$(d(S),d(S))$ is
on $N(S)$, we have $\alpha = (1 + m)d(S)$. Since  $m > 1$ when $m \geq 1$ is not an integer, the 
maximum possible value of $y$ on the line $t_1 + mt_2 = (1 + m)d(S)$ for $t_1, t_2 \geq 0$ is 
${m + 1 \over m} d(S) < 2d(S)$. Thus the degree of $y^{\beta}R(y^c)$ is less than $2d(S)$, and
hence the degree of $R(y)$ is less than ${2d(S) \over c} \leq d(S)$. Hence the zeroes of $R(y)$ are
of order less than $d(S)$, implying the zeroes of 
of $S_e(1,y) = y^{\beta}R(y^c)$ other than $y = 0$ are of order less than $d(S)$. This contradicts our
assumption that $S_e(1,y)$ has a zero $r \neq 0$ of order $k \geq d(S)$ and we conclude that $m$ 
is an integer. 

Note that since $m$ is an integer so is $\alpha$. Consequently by $(7.2')$ if $m$ is even,
then $S_e(1,y) = \pm S_e(-1,y)$, while if $m$ is odd one
has $S_e(1,y) = \pm S_e(-1,-y)$. Hence in either case both $S_e(1,y)$ and $S_e(-1,y)$ have a zero of 
order $k$ not at the origin. This completes the proof of Lemma 7.2.

The next lemma is the crux of this section. To set it up, suppose $S(x,y)$ is not in 
superadapted coordinates and the bisectrix intersects the interior of an edge $e$. Then
since $S_e(1,y)$ or $S_e(-1,y)$ has a zero $r \neq 0$ of order $k \geq d(S)$, $S_e(x,y)$ has zeroes of 
order $k$ at any point on a curve of the form $y = rx^m$. Hence $S_e(x,1)$ or $S_e(x,-1)$ has a zero of 
order $k$ away from the origin. Thus we may switch the roles
of the $x$ and $y$ axes if we want and assume $e$ has equation $t_1 + mt_2 = \alpha$ for $m \geq 1$; by
Lemma 7.2 $m$ is an integer and $S_e(1,y)$ has a zero $r \neq 0$ of order $k$.

\noindent {\bf Lemma 7.3.} Suppose $S(x,y)$ is not in superadapted coordinates and the bisectrix 
intersects the interior of an edge $e$. As described above, switching the $x$ and $y$ axes if
necessary, write the equation of $e$ as $t_1 + mt_2 = \alpha$ for an integer
$m \geq 1$ and assume $S_e(1,y)$ has a zero $r \neq 0$ of order $k \geq d(S)$. 
Then there is a coordinate change of the form $(x,y) \rightarrow (x,y + a(x))$ such that $a(x)$ is
smooth with $a(0) = 0$, after which one is either case 1 or 3 
of superadapted coordinates, or the following more general version of case 2:

\noindent {\bf Case 2'.} The bisectrix intersects $N(S)$ at a vertex $(d,d)$.

\noindent {\bf Proof.} Let $Q(y) = S_e(1,y)$, and let $(p,q)$ denote the upper vertex of
the edge $e$; necessarily $q > d(S)$. 
We will find a smooth function $a(x)$ such that  
$S'(x,y) = S(x,y + a(x))$ is in one of the following two mutually exclusive categories.

\noindent {\bf Category 1:} $S'(x,y)$ is either in case 1, case 2', or case 3.

\noindent {\bf Category 2:} The bisectrix intersects the interior of an edge $e'$ of $N(S')$
with equation $t_1 + m't_2 = \alpha '$, $m' > m \geq 1$, such that the upper vertex $(p',q')$ of $e'$ 
satisfies $q' < q$ and such that $S'(x,y)$ is not in case 1. (In particular by Lemma 7.2 $S'_{e'}
(1,y)$ has a zero of order $\geq d(S')$). 

\noindent Lemma 7.3 will then follow; there can be at most $q$ iterations of category 2.

We first consider the case where $k < q$. 
The function $Q(y + r)$ has a root at $y = 0$ of order $k$. We choose $a(x) = rx^m$ and define 
$S'(x,y) = S(x, y + a(x)) = S(x,y + rx^m)$. Note that $t_1 + mt_2 = \alpha$ is a supporting line of 
$N(S')$ as it was for $N(S)$, and that there is  
an edge $E$ of $N(S')$ on this line whose upper vertex is $(p,q)$. Observe that $S'_E(x,y) = S_e(x,y + 
rx^m) = x^{\alpha}Q({y \over x^m} + r)$. Since $Q$ has a zero of order $k$ at $r$, the lowest power of 
$y$ appearing in $S'_E(x,y)$ is $y^k$ and 
therefore $E$'s lower vertex is at a point $(j,k)$ for some $j$. Since both vertices of $E$ have 
$y$-coordinates
at least $d(S)$, they are both in the portion of the line $t_1 + mt_2 = \alpha$ on or above $(d(S),d(S))$. Thus
the edge $E$ lies wholly on or above the bisectrix. If the bisectrix intersects
$N(S')$ at a vertex or inside the horizontal or vertical rays, one is in Category 1. Otherwise, 
it must intersect $N(S')$ in the 
interior of an edge $e'$ whose upper vertex is either $(j,k)$ or a lower vertex. And because
$t_1 + mt_2 = \alpha$ is a supporting line for $N(S')$ and $e'$ lies below $E$, $e'$ will have equation 
$t_1 + m't_2 = \alpha '$ for some $m' > m \geq 1$. Thus we are either in case 1 superadapted coordinates 
(which is in Category 1) or we are in Category 2. Hence when $k < q$, $S'(x,y)$ is in 
either Category 1 or 2 and we are done.

It remains to consider the situation where $r$ is a zero of $Q(y)$ of order $q$. In this case we have
$Q(y) = c(y - r)^q$ for some $c$. For a large integer $n$ we expand $S(x,y)$ as
$$S(x,y) = cx^{\alpha}({y \over x^m} - r)^q + T_n(x,y) + E_n(x,y) \eqno (7.3)$$
Here the polynomial $T_n(x,y)$ are the terms of $S$'s Taylor expansion with exponents less
than $n$. For all $0 \leq \beta, \gamma < n$ one has
$$ |{\partial^{\beta + \gamma}E_n \over \partial x^{\beta} \partial y^{\gamma}}(x,y)| < C(|x|^{n - \beta} + 
|y|^{n - \gamma}) \eqno (7.4)$$
Note that
$$S(x,x^my) =cx^{\alpha}(y - r)^q + x^{\alpha + 1}T_n'(x,y) + E_n(x,x^my) \eqno
(7.5)$$
Here $T_n'(x,y)$ is also a polynomial. We define
$s(x,y) = {S(x,x^my) \over x^{\alpha}}$, so that
$$s(x,y) = c(y - r)^q + xT_n'(x,y) + x^{-\alpha}E_n(x,x^my) \eqno (7.6)$$
We claim that the function $s(x,y)$ is smooth on a neighborhood of $(0,r)$. Off the $y$-axis smoothness
holds because $S(x,y)$ is smooth. One can show that a given derivative of $s(x,y)$ exists when
$x = 0$ and equals that of $c(y - r)^q + xT_n'(x,y)$ for large enough $n$ by
examining the difference quotient of a one-lower order derivative of $(7.6)$, inductively assuming 
this lower-order derivative exists and has the right value when $x = 0$. Equation
$(7.4)$ ensures that the difference quotient of the lower derivative of $x^{-\alpha}E_n(x,x^my)$ tends 
to zero as $x$ goes to zero.
We conclude that $s(x,y)$ is smooth on a neighborhood of $(0,r)$. 

We next use the smooth implicit function 
theorem on ${\partial^{q-1}s \over \partial y^{q-1}}$ and find a smooth function $k(x)$ defined
in a neighborhood of $x = 0$ such that $k(0) = r$ and
${\partial^{q-1}s \over \partial y^{q-1}} (x,k(x)) = 0$.
Transferring this back to $S(x,y)$ we have
$${\partial^{q-1}S \over \partial y^{q-1}} (x,x^mk(x)) = 0 \eqno (7.7)$$
Thus if we let $a(x) = x^mk(x)$ and $S'(x,y) = S(x,y + x^mk(x))$, for all $x$ we consequently have
$${\partial^{q-1}S' \over \partial y^{q-1}} (x,0) = 0 \eqno (7.8)$$
Thus for every $a$ the Taylor series coefficient $S'_{a\,q-1}$ is zero.

Next, since $t_1 + mt_2 = 
\alpha$ is a supporting line for $N(S)$, this line is also a supporting
line for $N(S')$ and intersects $N(S')$ at the single vertex $(p,q)$. If $S'(x,y)$ is in Category 1 we
have nothing to prove, so we may assume we are not in Category 1. Let $e'$ denote the
edge of $N(S')$ intersecting the bisectrix and denote its equation by $t_1 + m't_2 = \alpha '$. 
Since $e'$ lies within the set $t_1 + mt_2 \geq \alpha$
and is no higher than the vertex $(p,q)$ of $N(S')$ that is on the supporting line $t_1 + mt_2 = \alpha$, 
we have $m' > m \geq 1$. If the upper vertex $(p',q')$ of $e'$ satisfies 
$q' < q$, one is in Category 2 and we are done. So we assume this upper vertex is $(p,q)$ itself.

If $S'_{e'}(1,y)$ has a real zero $r' \neq 0$ of order 
$k < q$, one is in the situation above $(7.3)$; there is a smooth $b(x)$ such that 
$S'(x,y + b(x)) = S(x,y + a(x) + b(x))$ is in Category 1 or 2 as needed. The only other 
possibility is that $S'_{e'}(1,y)$ has a single zero $r' \neq 0$ of order $q$. But 
this cannot happen. For this would imply $S'_{e'}(x,y) = c'x^{\alpha '}({y \over x^{m'}}- r')^q$ 
has a nonvanishing $y^{q-1}$ term. Consequently, for some $a$ the Taylor series
coefficient $S'_{a\,q-1}$ would be nonzero, contradicting $(7.8)$. Thus the case where 
$S'_{e'}(1,y)$ has a single zero of order $q$ does not occur, and we are done with the proof of Lemma 7.3.

The final step of the proof of the existence of superadapted coordinates is the following.

\noindent {\bf Lemma 7.4.} Suppose one is in case 2'; that is, the bisectrix intersects $N(S)$ at a 
vertex $(d,d)$. Then there exists a smooth coordinate change fixing the origin after which one is in
case 2 of superadapted coordinates.

\noindent {\bf Proof.} Suppose the bisectrix intersects $N(S)$ at a vertex $(d,d)$ but $S(x,y)$ is
not in superadapted coordinates. Then $(d,d)$ is on an edge $e$ of $N(S)$ such that $S_e(1,y)$ or 
$S_e(-1,y)$ has a 
zero $r \neq 0$ of order $d$ or greater. By Lemma 7.1, $S_e(x,y)$ is of the form $cx^{\alpha}({y \over x^m} - r)^d$ 
or $cy^{\alpha}({x \over y^m} - r)^d$ for positive integers $\alpha, m$ 
and some nonzero $c,r$, the first corresponding to the case where $e$ lies below $(d,d)$ and the second
corresponding to where it lies above $(d,d)$. Switching axes if necessary, assume that $S_e(x,y) = 
cx^{\alpha}({y \over x^m} - r)^d$. One can argue as in $(7.3)-(7.8)$ to
obtain a function of the form $T(x,y) = S(x, y + a(x))$ such that the bisectrix intersects 
$N(T)$ at the vertex $(d,d)$, but such that as in the last paragraph of the proof of Lemma 7.3 if $(d,d)$ 
is the upper vertex of an edge $e'$ of $N(T)$ then $T_{e'}(1,y)$ does not have a zero of order $d$ (or
greater). We 
also must have that $T_{e'}(-1,y)$ has no zero of order $d$ or greater; for if it did by Lemma 7.1 
we could write $T_{e'}(x,y)$ in the form $cx^{\alpha}({y \over x^m} - r)^d$, which would 
imply $T_{e'}(1,y)$ also has such a zero, a contradiction. 

There still remains the possibility that $(d,d)$ is the lower vertex of a compact edge $f$ of $N(T)$ 
such that $T_f(1,y)$ has a zero of order order $d$ or greater. By Lemma 7.1, if this happens
$T_f(x,y)$ is of the form $cy^a({x \over y^g} - r)^d$ where $g$ is an integer.
Again using the argument from $(7.3)$ onwards, this time reversing the roles of the $x$ and $y$ variables,
there is a smooth coordinate change of the form $\beta:(x,y) \rightarrow 
(x - y^gh(x,y),y)$ such that if one denotes $T \circ \beta$ by $S'$, then if $(d,d)$ is the 
lower vertex of an edge $f'$ of $N(S')$ 
then $S_{f'}'(x,1)$ and $S_{f'}'(x,-1)$ do not have any zeroes of order $d$ or greater other than $x = 0$. 
This means $S_{f'}'(1,y)$ and $S_{f'}'(-1,y)$ also have no such zero. For if one of the functions did, 
$S_{f'}'(x,y)$ would have zeroes of order $d$ along some curve $y = 
sx^n$, which would imply either $S_{f'}'(x,1)$ or $S_{f'}'(x,-1)$ had zeroes of order $d$ or greater 
away from the origin, a contradiction. 

Furthermore, the slope of $e'$ is of the form
${-{1 \over m'}}$ for $m'$ a positive integer and $g > {1 \over m'}$. As a result, $e'$ is an edge of
$N(S')$ containing $(d,d)$ and the coordinate change $\beta$ did not change any of the terms of $T_{e'}
(x,y)$. Thus $S_{e'}'(x,y) = T_{e'}(x,y)$ and all zeroes of $S_{e'}'(1,y)$ and $S_{e'}'(-1,y)$ other than
$y = 0$ have order less than $d$. Hence we are in superadapted coordinates and the proof of Theorem 7.4 
is complete.

\noindent {\bf 8. References.}

\noindent [ArKaCu] G.I. Arkhipov, A.A. Karacuba and V.N. Cubarikov, {\it Trigonometric Integrals},
Izv. Akad. Nauk SSSR Ser. Mat., {\bf 43}, no.5, (1979) 971-1003. \parskip = 3pt\baselineskip = 3pt

\noindent [AGV] V. Arnold, S Gusein-Zade, A Varchenko, {\it Singularities of differentiable maps},
Volume II, Birkhauser, Basel, 1988.

\noindent [C] M. Christ, {\it Hilbert transforms along curves. I. Nilpotent groups}, Annals of 
Mathematics (2) {\bf 122} (1985), no.3, 575-596.

\noindent [CaCW] A. Carbery, M. Christ, J Wright, {\it Multidimensional van der Corput and sublevel
set estimates}, J. Amer. Math. Soc. 12 (1999), no. 4, 981--1015. 

\noindent [DNS] J. Denef,  J. Nicaise, P. Sargos, {\it Oscillating integrals and Newton polyhedra},
J. Anal. Math. {\bf 95} (2005), 147-172.

\noindent [DS] J. Denef, P, Sargos, {\it Polyedre de Newton et distribution $f\sp s\sb +$. II.}
Math. Ann. {\bf 293} (1992), no. 2, 193-211.

\noindent [F] M.V. Fedoryuk, {\it The saddle-point method}, Nauka, Moscow, 1977.

\noindent [G1] M. Greenblatt, {\it Resolution of singularities, asymptotic expansions of
integrals over sublevel sets, and applications}, submitted.

\noindent [G2] M. Greenblatt, {\it Sharp $L\sp 2$ estimates for one-dimensional 
oscillatory integral operators with $C\sp \infty$ phase.} Amer. J. Math. 
{\bf 127} (2005), no. 3, 659-695.

\noindent [IM] I. Ikromov, D. M\"uller, {\it On adapted coordinate systems}, preprint.

\noindent [IKM] I. Ikromov, M. Kempe, and D. M\"uller, {\it Sharp $L^p$ estimates for maximal operators
associated to hypersurfaces in $\R^3$ for $p > 2$}, preprint.

\noindent [R] V. Rychkov, {\it Sharp $L^2$ bounds for oscillatory
integral operators with $C^\infty$ phases}, Math. Zeitschrift, {\bf 236}
(2001) 461-489.

\noindent [PS] D. H. Phong, E. M. Stein, {\it The Newton polyhedron and
oscillatory integral operators}, Acta Mathematica {\bf 179} (1997), 107-152.

\noindent [PSSt] D. H. Phong, E. M. Stein, J. Sturm, {\it On the growth and 
stability of real-analytic functions}, Amer. J. Math. {\bf 121} (1999), no. 3, 519-554.

\noindent [Sh] H. Schulz, {\it Convex hypersurfaces of finite type and Their Fourier Transforms}, Indiana 
Univ. Math. Journal, {\bf 40} no. 4 (1991), 1267-1275.

\noindent [S] E. Stein, {\it Harmonic analysis; real-variable methods, orthogonality, and oscillatory
integrals}, Princeton Mathematics Series Vol. 43, Princeton University Press, Princeton, NJ, 1993.

\noindent [V] A. N. Varchenko, {\it Newton polyhedra and estimates of
oscillatory integrals}, Functional Anal. Appl. {\bf 18} (1976), no. 3, 
175-196.

\noindent [Vi] I. M. Vinogradov, {\it The method of trigonometrical sums in number theory},
Trav. Inst. Math. Stekloff  {\bf 23},  (1947). 109 pp. 

\line{}
\line{}

\noindent Department of Mathematics, Statistics, and Computer Science \hfill \break
\noindent University of Illinois at Chicago \hfill \break
\noindent 322 Science and Engineering Offices \hfill \break
\noindent 851 S. Morgan Street \hfill \break
\noindent Chicago, IL 60607-7045 \hfill \break
\end